%% file: d0-cnfd.tex
\documentclass[11pt]{article}

\usepackage{amsmath, amsthm}
\usepackage{amsmath, amsfonts}
\usepackage{amsmath, amssymb}
\usepackage{amsmath}
\usepackage[all]{xy}

\newtheorem{definition-lemma}[theorem]{Definition/Lemma}
\newtheorem{definition-example}[theorem]{Definition/Example}
\newtheorem{definition-explanation}[theorem]{Definition/Explanation}
\newtheorem{explanation-definition}[theorem]{Explanation/Definition}

\newtheorem{lemma-definition}[theorem]{Lemma/Definition}

\numberwithin{equation}{subsection}

\newtheorem{stheorem}{Theorem}[section]
\newtheorem{sdefinition}[stheorem]{Definition}
\newtheorem{sdefinition-lemma}[stheorem]{Definition/Lemma}
\newtheorem{sdefinition-example}[stheorem]{Definition/Example}
\newtheorem{sdefinition-explanation}[stheorem]{Definition/Explanation}
\newtheorem{sexplanation-definition}[stheorem]{Explanation/Definition}
\newtheorem{slemma}[stheorem]{Lemma}
\newtheorem{slemma-definition}[stheorem]{Lemma/Definition}
\newtheorem{sproposition}[stheorem]{Proposition}

\newtheorem{sremark}[stheorem]{\it Remark}

\newtheorem{sexample}[stheorem]{Example}

\input{epsfig.sty}
\input{newcmd.tex}

\begin{document}

\enlargethispage{23cm}

\begin{titlepage}

$ $

\vspace{-1cm}
% \vspace{-1.5cm} % Re: -1.5cm for PC; -2.5cm for UT-Math-system

\noindent\hspace{-1cm}
\parbox{6cm}{\small May 2009}\
   \hspace{8cm}\
   \parbox[t]{5cm}{yymm.nnnn [math.AG]\\ D(4): D$0$, conifold.}

\vspace{2cm}

%title
\centerline{\large\bf
 Azumaya structure on D-branes}
\vspace{1ex}
\centerline{\large\bf
 and deformations and resolutions of a conifold revisited:}
\vspace{1ex}
\centerline{\large\bf
 Klebanov-Strassler-Witten vs.\ Polchinski-Grothendieck}
% end-title

\bigskip

\vspace{3em}

%authors-'n-addresses
\centerline{\large
  Chien-Hao Liu
  \hspace{1ex} and \hspace{1ex}
  Shing-Tung Yau
}

\vspace{6em}

%abstract%
\begin{quotation}
\centerline{\bf Abstract}

\vspace{0.3cm}

\baselineskip 12pt  %13pt for [12pt] style
{\small
  In this sequel to [L-Y1], [L-L-S-Y], and [L-Y2]
  (respectively arXiv:0709.1515 [math.AG],
        arXiv:0809.2121 [math.AG], and
        arXiv:0901.0342 [math.AG]),
  we study a D-brane probe on a conifold
   from the viewpoint of the Azumaya structure on D-branes
   and toric geometry.
  The details of how deformations and resolutions of
   the standard toric conifold $Y$ can be obtained
   via morphisms from Azumaya points are given.
  This should be compared with the quantum-field-theoretic/D-brany
   picture of deformations and resolutions of a conifold
   via a D-brane probe sitting at the conifold singularity
   in the work of Klebanov and Witten [K-W] (arXiv:hep-th/9807080)
   and Klebanov and Strasser [K-S] (arXiv:hep-th/0007191).
  A comparison with resolutions via noncommutative desingularizations
   is given in the end.
} % endsmall
\end{quotation}

%\bigskip
\vspace{12em}

\baselineskip 12pt
{\footnotesize
} %end-footnotesize

\bigskip

\noindent {\small MSC number 2000:
 14E15, 81T30; 14A22, 16G30, 81T75.
} % end-small

\bigskip

\baselineskip 10pt
% Re: 11pt for [11pt] style; 12pt for [12pt] style
{\scriptsize
\noindent{\bf Acknowledgements.}
 We thank
  Kwokwai Chan, Miranda Chih-Ning Cheng, Peng Zhang
   for discussions in the mirror symmetry seminar.
 C.-H.L.\ thanks in addition
  Cumrun Vafa
   for numerous illuminations on branes
   and answer-to-questions throughout his topic course in spring 2009;
  Alessandro Tomasiello
   for two long discussions after his substitute lectures for C.V.;
  Clay Cordova
   for the string discussion session;
  participants of the string reading seminar on AdS/CFT, organized by
   Chi-Ming Chang;
  Clark Barwick, Yum-Tong Siu for other topic courses;
  and Ling-Miao Chou for moral support.
 The project is supported by NSF grants DMS-9803347 and DMS-0074329.
} %endscriptsize

\end{titlepage}

\newpage
\begin{titlepage}
\baselineskip 14pt

$ $

\vspace{12em} % \vspace{4em}

\centerline{\it
 In memory of a young string theorist Ti-Ming Chiang,}
\vspace{.6ex}
\centerline{\it
 whose path I crossed accidentally and so briefly.$^{\dagger}$}

\vspace{30em}

\baselineskip 10pt
 {\scriptsize
\noindent$^{\dagger}${\it From C.-H.L.}\hspace{1em}
 During the years I was attending Prof.\ Candelas's group meetings,
  I learned more about Calabi-Yau manifolds and mirror symmetry and
  got very fascinated by the works from Brian Greene's group.
 Because of this, I felt particularly lucky
   knowing later that I was going to meet one of his students,
   Ti-Ming, - a young string theorist with a PhD from Cornell
    at his very early 20's -
  and perhaps to cooperate with him.
 Unfortunately that anticipated cooperation never happened.
  Ti-Ming had become unwell just before I resettled.
 Except the visits to him at the hospital and some chats
  when he showed up in the office,
 I didn't really get the opportunity to interact with him intellectually.
 Further afterwards I was informed of Ti-Ming's passing away.
 Like a shooting star he reveals his shining so briefly
  and then disappears.
 The current work is the last piece of Part 1 of
  the D-brane project.
 It is grouped with the earlier D(1), D(2), D(3)
  under the hidden collective title:
  ``Azumaya structure on D-branes and its tests".
 Here we address in particular a conifold from the viewpoint
  of a D-brane probe with an Azumaya structure.
 This is a theme Ti-Ming may have felt interested in as well,
   should he still work on string theory,
  since
   conifolds have play a role in understanding the duality web
    of Calabi-Yau threefolds - a theme Ti-Ming once worked on -
    and D-brane resolution of singularities is a theme
    Brian Greene's group once pursued vigorously.
 We thus dedicate this work to the memory of Ti-Ming.
} %endscriptsize

\end{titlepage}

%paper
\newpage
$ $

\vspace{-4em}  % Re: -4cm for PC; -6cm for UT-Math-system

%short heading
\centerline{\sc Deformations and Resolutions of a Conifold
                via a D-Brane Probe}

\vspace{2em}

\baselineskip 14pt  %Re: 14pt for [11pt] style
                    %Re: 15pt for [12pt] style.

\begin{flushleft}
{\Large\bf 0. Introduction and outline.}
\end{flushleft}
Conifolds,
i.e.\ Calabi-Yau threefolds with ordinary double-points,
 have been playing special roles
 at various stages of string theory.\footnote{Readers are
                    referred to, for example,
                     [C-dlO] (1989);
                     [Stro], [G-M-S], [C-G-G-K] (1995);
                     [G-V] (1998);
                     [Be], [C-F-I-K-V] (2001) and references therein
                     to get a glimpse of conifolds in string theory
                     around the decade 1990s.}
In this sequel to [L-Y1], [L-L-S-Y], and [L-Y2],
 we study a D-brane probe on a conifold
 from the viewpoint of Azumaya structure on D-branes
 and toric geometry.
This should be compared with the quantum-field-theoretic/D-brany
 picture of deformations and resolutions of a conifold
   in the work of Klebanov and Strasser [K-S] and
   Klebanov and Witten [K-W].

\bigskip

\begin{flushleft}
{\bf Effective-space-time-filling D$3$-brane at a conifold singularity.}
\end{flushleft}
In [K-W], Klebanov and Witten studied
 the $d=4$, $N=1$ superconformal field theory (SCFT)\footnote{There
                 will be a few standard physicists' conventional notations
                  in this highlight of the relevant part of
                 [K-W] and [K-S]:
                 $N$ that counts the {\it number of supersymmetries}
                  (susy) via the multiple number of minimal susy numbers
                  in each space-time dimension  vs.\
                 $N$ that appears in the {\it gauge group}
                  $U(N)$ or {\it SU}$\,(N)$ vs.\
                 $N$ that counts the {\it multiplicity} of
                  stacked D-branes.}
  on the D$3$-brane world-volume $X$
   ($\simeq {\Bbb R}^4$ topologically)
  that is embedded in the product space-time ${\Bbb M}^{3+1}\times Y$
  as ${\Bbb M}^{3+1}\times\{{\mathbf 0}\}$,\footnote{In
                     string-theorist's terminology, the D$3$-brane
                      is ``{\it sitting at the conifold singularity}".
                     We will also adopt this phrasing for convenience.
                     Note that in such a setting, the internal part
                      is a D$0$-brane on the conifold $Y$.
                     The latter is what we will study in this work. }
  and its supergravity dual
  - a compactification of $d=10$, type-IIB supergravity theory
    on AdS$^5\times (S^3\times S^2)$ -
  along the line of the AdS/CFT correspondence of Maldacena [Ma].
Here
 ${\Bbb M}^{3+1}$ is the $d=3+1$ Minkowski space-time,
 $Y$ is the conifold $\{z_1z_2-z_3z_4=0\}\subset {\Bbb C}^4$
  (with coordinates $(z_1, z_2,z_3, z_4)$),
 ${\mathbf 0}$ is the conifold singularity on $Y$, and
 AdS$^5$ is the $d=4+1$ anti-de Sitter space-time.

In the simplest case when there is a single D$3$-brane
  sitting at the conifold point of $Y$,
 the {\it classical moduli space} of the {\it supersymmetric vacua}
  of the associated $U(1)$ super-Yang-Mills theory coupled with matter
   on the D$3$-brane world-volume
  comes from the {\it $D$-term} of the {\it vector multiplet} and
  the coefficient $\zeta\in {\Bbb R}$ of
  the {\it Fayet-Iliopoulos term} in the Lagrangian.\footnote{$\zeta$
                    is part of the parameters to give local coordinates
                     of the {\it Wilson's theory-space} in the problem;
                    cf.~[L-Y2: Introduction] for brief words.
                    See also [W-B] for the standard SUSY jargon.}
By varying $\zeta$, one realizes the two small resolutions,
 $Y_+$ and $Y_-$, of $Y$ as the classical moduli space $Y_{\zeta}$
 of the above $d=4$ SCFT.\footnote{See
                              also [Wi] and [D-M] for details
                              of such a construction.}
A flop $\xymatrix{X_+ \ar @{-->}[r] & Y_-}$ happens
 when $Y_{\zeta}$ crosses over $\zeta=0$.

To describe the physics for $N$-many parallel D$3$-branes
 sitting at the conifold singularity,
Klebanov and Witten proposed to enlarge the gauge group for
 the super-Yang-Mills theory on the common world-volume of
 the stacked D$3$-brane to $U(N)\times U(N)$
  (rather than the naive $U(N)$)  {\it and}
 introduce a {\it superpotential} $W$ for the chiral multiplets.
The classical moduli space of the theory comes from a system
 with equations of the type above (i.e.\ {\it D-term equations})
  and
 equations from the superpotential term $W$
  (i.e.\ {\it F-term equations}).
In particular, the $N$-fold symmetric product $\Sym^nY$ of $Y$
  can be realized as the classical moduli space of
  the $d=4$ SCFT on the D$3$-brane world-volume with $\zeta=0$.

In [K-S], Klebanov-Strassler studied further $d=4$, $N=1$
 supersymmetric quantum field theory (SQFT)
 on the D$3$-brane world-volume that arises from
 a D$3$-brane configuration
  with both $N$-many above {\it full}/free D$3$-branes and
       $M$-many new {\it fractional}/trapped D$3$-branes\footnote{See
                                      [G-K] and references therein
                                       for the detail of such
                                       {\it fractional D-branes}.}
  sitting at the conifold singularity ${\mathbf 0}$ of $Y$.
For infrared physics, the theory now has the gauge group
 $\SU(N+M)\times \SU(N)$.
It follows from the work of Affleck, Dine, and Seiberg
                      [A-D-S]\footnote{See
                           also [Arg: Chapter~3] and [Te: Chapter~9].}
 that an {\it additional term} to the previous superpotential $W$
  is now {\it dynamically generated}.
This deforms the classical moduli space of SUSY vacua of
 the $d=4$ SQFT on the D$3$-brane world-volume.
In the simplest case when $N=M=1$,
 this enforces a deformation of the classical moduli space
 from a conifold to a deformed conifold $Y^{\prime}$
 ($\simeq T^{\ast}S^3$ topologically).
Cf.\ {\sc Figure}~0-1.

\begin{figure}[htbp]
 \epsfig{figure=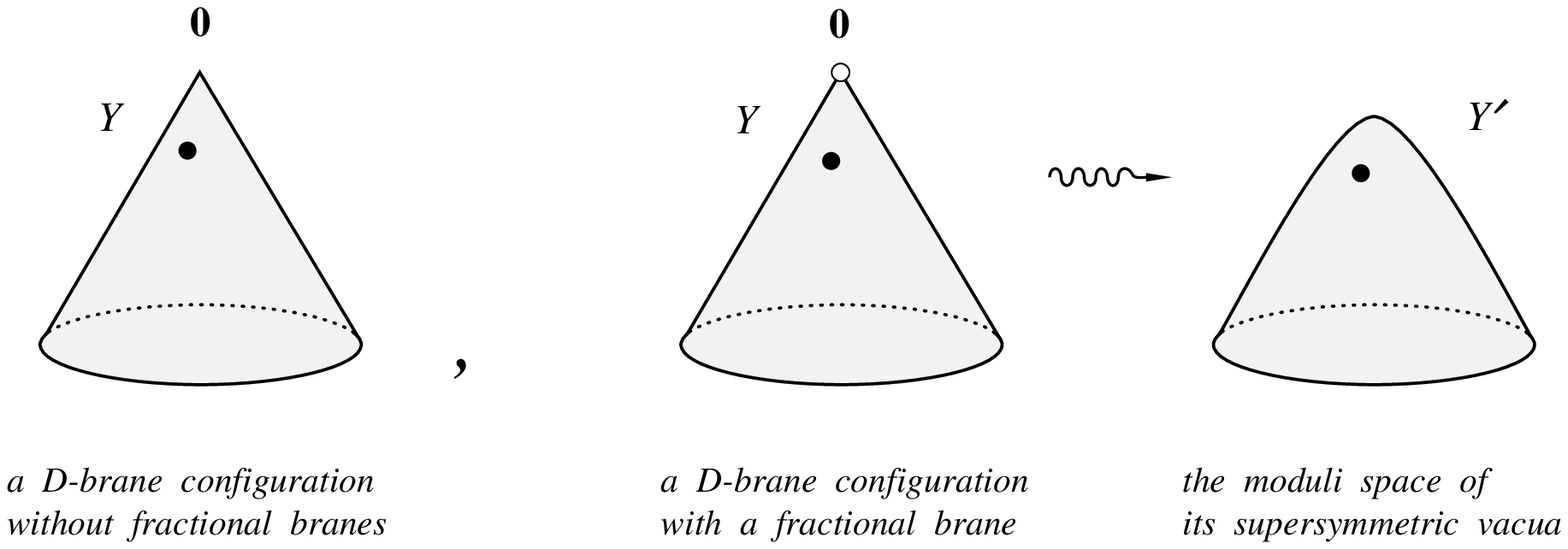,width=16cm}
 \centerline{\parbox{13cm}{\small\baselineskip 12pt
  {\sc Figure}~0-1. (Cf.\ [Stra: Figures 25, 26, 27].)
  When a fractional/trapped D$3$-brane sits
   at the conifold singularity ${\mathbf 0}\in Y$,
  the full/free D$3$-brane ``sees" a smooth deformed conifold
   $Y^{\prime}$ ($\simeq T^{\ast}S^3$ topologically)
  as its classical vacua manifold.
 I.e., in very low energy for this situation
  the free D$3$-brane ``feels"
  as if it lives on $Y^{\prime}$ instead of $Y$!
 In the figure, a full D$3$-brane is indicated by $\bullet$ while
  a fractional D$3$-brane by $\circ$.
 }}
\end{figure}

While giving only a highlight of key points in [K-W] and [K-S]
 that are most relevant to us,
we should remark that,
 in addition to further quantum-field-theoretical issues
  on the gauge theory side,
 there is also a gravity side of the story that was studied
  in [K-W] and [K-S].\footnote{See
          [A-G-M-O-O] and [Stra] for a review with more emphasis on
           respectively the gravity and the gauge theory side
           in the correspondence;
          e.g.\ [G-K], [K-N] for developments between [K-W] and [K-S];
          and e.g.\ [D-K-S] for a more recent study.}

\bigskip

\begin{flushleft}
{\bf Azumaya structure on D-branes and its tests.}
\end{flushleft}
In D(1) [L-Y1], D(2) [L-L-S-Y], D(3) [L-Y2] and the current work D(4),
 we illuminate the Azumaya geometry
 as a key feature of the geometry on D-brane world-volumes
 in the algebro-geometric category.
These four together center around
 the very remark of Polchinski:
 \begin{quote}
 ([Po: vol.~1, Sec.~8.7, p.~272]) \hspace{1em}
  ``{\it
  For the collective coordinate $X^{\mu}$, however, the meaning
   is mysterious: the collective coordinates for the embedding of
   $n$ D-branes in space-time are now enlarged to $n\times n$ matrices.
  This `noncommutative geometry' has proven to play a key role in
   the dynamics of D-branes, and there are conjectures that
   it is an important hint about the nature of space-time.}",
 \end{quote}
 which was taken as a guiding question as to what a D-brane is
 in this project, cf.\ [L-Y1: Sec.~2.2].
D(2), D(3), and the current D(4)
 are meant to give more explanations of the highlight [L-Y1: Sec.~4.5].
In this consecutive series of four, we learned that$\,$:

\bigskip

\noindent
{\bf Lesson 0.1 [Azumaya structure on D-branes].} {\it
 This ``enhancement to $n\times n$ matrices" Polchinski alluded to
  says even more fundamentally the nature of D-branes themselves,
  i.e.\ the Azumaya structure thereupon.
 This structure gives them
  the power to detect the nature of space-time.
 We also learned that
 Azumaya structures on D-branes and morphisms therefrom
  can be used to reproduce/explain several stringy/brany phenomena
  of stringy or quantum-field-theoretical origin
  that are very surprising/mysterious at a first mathematical glance.
} % end-lesson

\bigskip

\noindent
This is a basic test to ourselves to believe that
 Azumaya structures play a special role in
 understanding/desccribing D-branes in string theory.
Having said this,
 we should however mention that
 D-brane remains a very complicated object and
 the Azumaya structure addressed here is only a part of it.
Further issues are investigated in separate works.

\bigskip

\noindent
{\bf Convention.}
 Standard notations, terminology, operations, facts in
  (1) physics aspects of strings and D-branes;
  (2) algebraic geometry;
  (3) toric geometry
  can be found respectively in
  (1) [Po], [Jo];
  (2) [Ha];
  (3) [Fu].
 \begin{itemize}
  \item[$\cdot$]
   {\it Noncommutative algebraic geometry} is a very technical topic.
   For the current work,
     [Art] of Artin,
     [K-R] of Kontsevich and Rosenberg,
     and [leB1] of Le Bruyn
    are particularly relevant.
   See [L-Y1: References] for more references.
 \end{itemize}

\bigskip

\begin{flushleft}
{\bf Outline.}
\end{flushleft}
{\small
\baselineskip 12pt  %13pt
\begin{itemize}
 \item[0.]
  Introduction.
  \vspace{-.6ex}
  \begin{itemize}
   \item[$\cdot$]
    Effective-space-time-filling D$3$-brane at a conifold singularity

   \item[$\cdot$]
    Azumaya structure on D-branes and its tests
  \end{itemize}

 \item[1.]
  D-branes in an affine noncommutative space.
  \vspace{-.6ex}
  \begin{itemize}
   \item[$\cdot$]
    Affine noncommutative spaces and their morphisms

   \item[$\cdot$]
    D-branes in an affine noncommutative space
    \`{a} la Polchinski-Grothendieck Ansatz
  \end{itemize}

 \item[2.]
  Deformations of a conifold via an Azumaya probe.
  \vspace{-.6ex}
  \begin{itemize}
   \item[$\cdot$]
    a toric setup for the standard local conifold

   \item[$\cdot$]
    an Azumaya probe to a noncommutative space and
    its commutative descent

   \item[$\cdot$]
    deformations of the conifold via an Azumaya probe:\\
    descent of noncommutative superficially-infinitesimal deformations.

   \item[$\cdot$]
    deformations of the conifold via an Azumaya probe: details
  \end{itemize}

 \item[3.]
  Resolutions of a conifold via an Azumaya probe.
  \vspace{-.6ex}
  \begin{itemize}
   \item[$\cdot$]
    D-brane probe resolutions of a conifold via the Azumaya structure

   \item[$\cdot$]
    an explicit construction of
     $\widetilde{Y}^{\prime}$, $Y_+^{\prime}$, and $Y_-^{\prime}$

   \item[$\cdot$]
    a comparison with resolutions via noncommutative desingularizations
  \end{itemize}
\end{itemize}
} %endsmall

\newpage
\section{D-branes in an affine noncommutative space.}

We recall definitions and notions in [L-Y1]
 that are needed for the current work.
Readers are referred to ibidem for more details and references.
See also [L-L-S-Y] and [L-Y2] for further explanations and examples.

\bigskip

\begin{flushleft}
{\bf Affine noncommutative spaces and their morphisms.}
\end{flushleft}
An {\it affine noncommutative space} over ${\Bbb C}$
 is meant to be a ``space" $\Space R$
  that is associated to an associative unital ${\Bbb C}$-algebra $R$.
In general, it can be tricky to truly realize $\Space R$
 as a set of points with a topology in a natural/functorial way.
However,``geometric" notions can still be pursued
  - despite not knowing what $\Space R$ really is -
 via imposing the fundamental geometry/algebra ansatz:
 \begin{itemize}
  \item[$\cdot$]
  {\it $[\,$geometry$\,=\,$algebra$\,]$}\hspace{1ex}
   The correspondence $\,R\leftrightarrow\Space R\,$
    gives a contravariant equivalence
    between the category $\AlgCategory_{\Bbb C}$
      of associative unital ${\Bbb C}$-algebras
     and the category $\AffineSpace_{\Bbb C}$
      of ``affine noncommutative spaces" over ${\Bbb C}$.
 \end{itemize}
For example,

\smallskip

\begin{sdefinition}
{\bf [smooth affine noncommutative space].}
{\rm ([C-Q: Sec.~3], [K-R: Sec.~1.1.4].)} {\rm
 An affine noncommutative space $\Space R$ over ${\Bbb C}$
  is said to be {\it smooth}
  if the associative unital ${\Bbb C}$-algebra $R$
   is finitely generated and
   satisfies the following property:
   \begin{itemize}
    \item[$\cdot$]
    ({\it lifting property for nilpotent extensions})\hspace{1em}
    for any ${\Bbb C}$-algebra $S$,
        two-sided nilpotent ideal $I\subset R$
         (i.e.\ $I=BIB$ and $I^n=0$ for $n>>0$), and
        ${\Bbb C}$-algebra homomorphism $h:R\rightarrow B/I$,
     there exists an ${\Bbb C}$-algebra homomorphism
     $\widetilde{h}:R\rightarrow S$ such that the diagram
     $$
     \xymatrix{
                                               && B \ar[d]\\
        R \ar[rr]_h \ar[rru]^{\widetilde{h}}   && B/I
     }
     $$
     commutes. Here $B\rightarrow B/I$ is the quotient map.
   \end{itemize}

}\end{sdefinition}

\smallskip

The following two classes of smooth affine noncommutative spaces
 are used in this work.

\smallskip

\begin{sexample}
{\bf [noncommutative affine space].}
{\rm ([K-R: Sec.~2: Example (E1)].)} {\rm
 The {\it noncommutative affine $n$-space}
  $\NA^n
   :=\Space ({\Bbb C}\langle\,\xi_1\,,\,\cdots\,,\,\xi_n\, \rangle)$
  {\it over ${\Bbb C}$} is smooth.
 Here ${\Bbb C}\langle\,\xi_1\,,\,\cdots\,,\,\xi_n\, \rangle$
  is the associative unital ${\Bbb C}$-algebra
  freely generated by the elements in the set
  $\{\,\xi_1\,,\,\cdots\,,\,\xi_n\,\}$.
}\end{sexample}

\begin{sexample}
{\bf [Azumaya-type noncommutative space].}
{\rm ([C-Q: Sec.~5 and Proposition~6.2],
      [K-R: Sec.~1.2, Examples (E2) and (C4)].)} {\rm
 Let $M_r(R)$ be the ${\Bbb C}$-algebra of $r\times r$-matrices
  with entries in a commutative regular ${\Bbb C}$-algebra $R$.
 Then the {\it Azumaya-type noncommutative space} $\Space M_r(R)$
  is smooth (over ${\Bbb C}$).
 Furthermore, it is also smooth over $\Spec R$.
}\end{sexample}

\smallskip

As a consequence of the Geometry/Algebra Ansatz,
a {\it morphism} $\varphi:X=\Space R \rightarrow Y=\Space S$
  is defined contravariantly
 to be a ${\Bbb C}$-algebra homomorphism
  $\varphi^{\sharp}: S\rightarrow R$.
The {\it image}, denoted $\Image\varphi$ or $\varphi(X)$,
  of $X$ under $\varphi$
 is defined to be $\Space(S/\Ker\varphi^{\sharp})$.
The latter is canonically included in $Y$
 via the morphism $\iota: \varphi(X) \hookrightarrow Y$
 defined by the ${\Bbb C}$-algebra quotient-homomorphism
 $\iota^{\sharp}: S\rightarrow S/\Ker\varphi^{\sharp}\,$.
This extends what is done in Grothendieck's theory
 of (commutative) schemes.
The benefit of thinking a morphism between affine noncommutative
 spaces this way is actually {\it two} folds:
 \begin{itemize}
  \item[(1)]
   {\it As a functor of point}$\,$:
   The space $X=\Space R$ defines a functor
    $$
     \begin{array}{ccccc}
      h_X  & : & \AffineSpace_{\Bbb C}  & \longrightarrow
           & \SetCategory^{\circ} \\[.6ex]
      && Y & \longmapsto  & \Mor(Y, X)\,;
     \end{array}
    $$
   i.e.\ a functor
    $$
     \begin{array}{ccccc}
      h_R  & : & \AlgCategory_{\Bbb C}  & \longrightarrow
           & \SetCategory \\[.6ex]
      && S & \longmapsto  & \Hom(R,S)\,.
     \end{array}
    $$
   Here $\SetCategory$ is the category of sets,
    $\SetCategory^{\circ}$ its opposite category,
    and $\Hom(R,S)$ is the set of ${\Bbb C}$-algebra-homomorphisms.

  \item[(2)]
   {\it As a probe}$\,:$
   $X=\Space R$ defines another functor
    $$
     \begin{array}{ccccc}
      g_X  & : & \AffineSpace_{\Bbb C}  & \longrightarrow
           & \SetCategory \\[.6ex]
      && Y & \longmapsto  & \Mor(X, Y)\,;
     \end{array}
    $$
   i.e.\ a functor
    $$
     \begin{array}{ccccc}
      g_R  & : & \AlgCategory_{\Bbb C}  & \longrightarrow
           & \SetCategory^{\circ} \\[.6ex]
      && S & \longmapsto  & \Hom(S,R)\,.
     \end{array}
    $$
 \end{itemize}
Aspect (1) is by now standard in algebraic geometry. It allows one
to define the various {\it local geometric properties}
 of a ``space" via algebra-homomorphisms;
for example, Definition~1.1.
 % Definition [smooth affine noncommutative space]
It suggests one to think of $X$ as a sheaf over $\AffineSpace_{\Bbb C}$.
Thus, after the notion of {\it coverings} and {\it gluings} is selected,
 it allows one to extend the notion of a noncommutative space to
 that of a ``{\it noncommutative stack}".
Aspect (2) is especially akin to our thought on D-branes.
It says, in particular, that the geometry of $X=\Space R$
 can be revealed through an ${\Bbb C}$-subalgebra of $R$.

\smallskip

\begin{sexample}
{\bf $[$Azumaya point$\,]$.} {\rm
 Consider the {\it Azumaya point of rank $r\,$}:
  $\Space M_r({\Bbb C})$.
 Its only two-sided prime ideal is $(0)$, the zero ideal.
 Thus, naively, one would expect $\Space M_r({\Bbb C})$
  to behave like a point with an Artin ${\Bbb C}$-algebra
  as its function ring.
 {\it However}, for example, from the ${\Bbb C}$-algebra monomorphism
  $\times^r{\Bbb C}\hookrightarrow M_r({\Bbb C})$
  with image the diagonal matrices in $M_r({\Bbb C})$,
  one sees that $\Space M_r({\Bbb C})$
    - which is topologically a one-point set
       if one adopts its interpretation as $\Spec M_r({\Bbb C})$ -
   can dominate $\amalg_r\Spec{\Bbb C}$
    - which is topologically a disjoint union of $r$-many points -.
 Furthermore, consider, for example, the morphism
  $\varphi:\Space M_r({\Bbb C})\rightarrow {\Bbb A}^1=\Spec {\Bbb C}[z]$
  defined by $\varphi^{\sharp}: {\Bbb C}[z]\rightarrow M_r({\Bbb C})$
   with $\varphi^{\sharp}(z) = m$ that is diagonalizable with
    $r$ distinct eigenvalues $\lambda_1,\,\cdots\,,\, \lambda_r$.
 Then $\Image\varphi$ is a collection of $r$-many ${\Bbb C}$-points
  on ${\Bbb A}^1$, located at $z=\lambda_1,\,\cdots\,,\,\lambda_r$
  respectively.
 In other words, the Azumaya noncommutativity cloud $M_r({\Bbb C})$
  over the seemingly one-point space $\Space M_r({\Bbb C})$
  can really ``split and condense" to
  a collection of concrete geometric points!
 Cf.~{\sc Figure}~1-1.
 See [L-Y1: Sec.~4.1] for more examples.
 Such phenomenon generalizes to Azumaya schemes;
  in particular, see [L-L-S-Y] for the case of Azumaya curves.
}\end{sexample}

\begin{figure}[htbp]
 \epsfig{figure=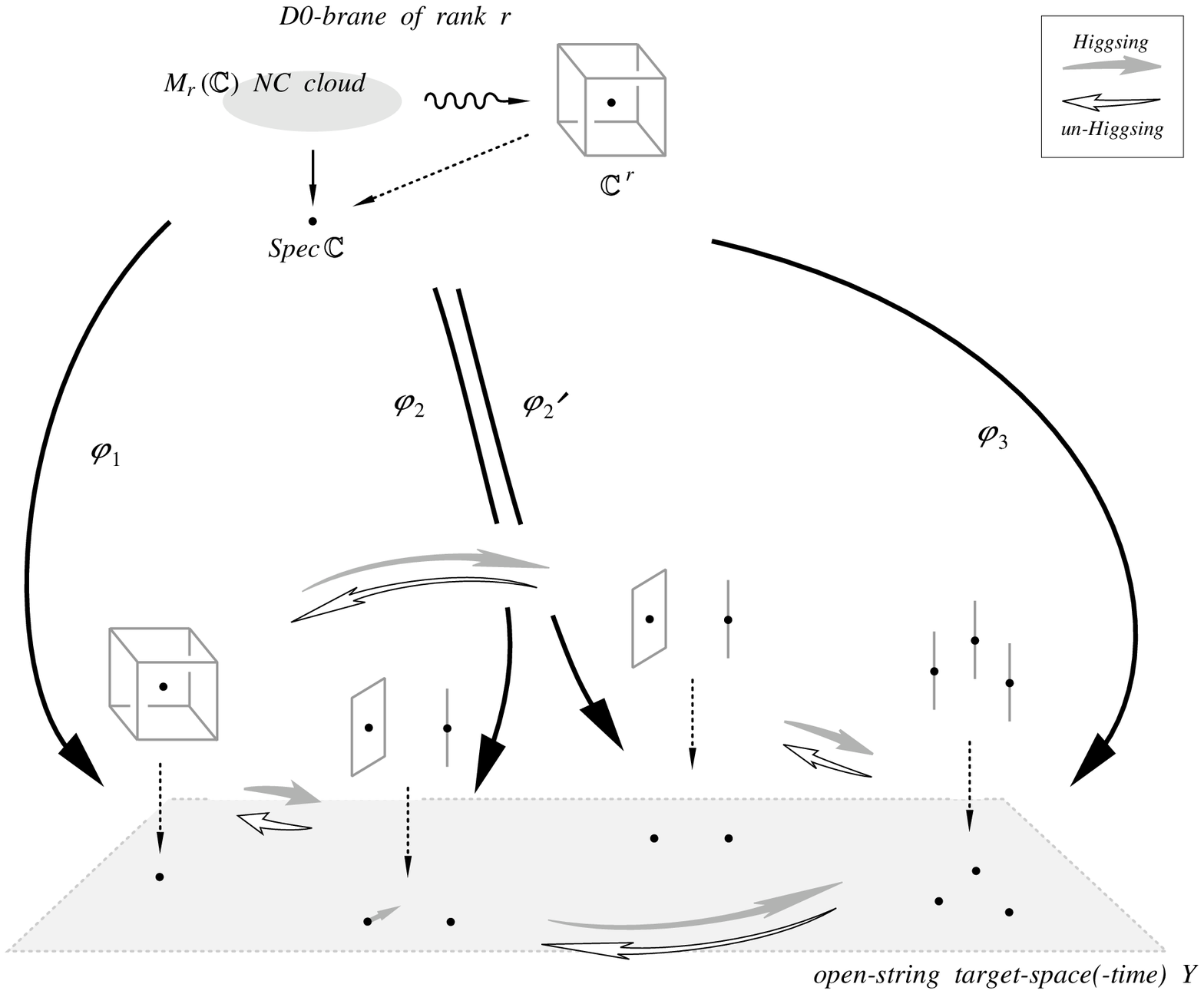,width=16cm}
 \centerline{\parbox{13cm}{\small\baselineskip 12pt
  {\sc Figure}~1-1. ([L-L-S-Y: {\sc Figure}~2-1-1].)
 Despite that {\it Space}$\,M_r({\Bbb C})$ may look only
   one-point-like,
  under morphisms
  the Azumaya ``noncommutative cloud" $M_r({\Bbb C})$
  over {\it Space}$\,M_r({\Bbb C})$ can ``split and condense"
  to various schemes with a rich geometry.
 The latter schemes can even have more than one component.
 The Higgsing/un-Higgsing behavior of the Chan-Paton module of
   D$0$-branes on $Y$ occurs
  due to the fact that
   when a morphism
    $\varphi:$ {\it Space}$\,M_r({\Bbb C}) \rightarrow Y$
    deforms,
   the corresponding push-forward $\varphi_{\ast}{\Bbb C}^r$
    of the fundamental module ${\Bbb C}^r$
    on {\it Space}$\,M_r({\Bbb C})$ can also change/deform.
 These features generalize to morphisms from Azumaya schemes to $Y$.
 Here, a module over a scheme is indicated by a dotted arrow
  $\xymatrix{ \ar @{.>}[r] &}$.
 }}
\end{figure}

\begin{sdefinition}
{\bf [surrogate associated to morphism].} {\rm
 Given $X=\Space R$,
  let $R^{\prime}\hookrightarrow R$ be a ${\Bbb C}$-subalgebra of $R$.
 Then, the space $X^{\prime}:=\Space R^{\prime}$
  is called a {\it surrogate} of $X$.
 By definition, there is a built-in dominant morphism
  $X\rightarrow X^{\prime}$,
  defined by the inclusion $R^{\prime}\hookrightarrow R$.
 Given a morphism $\varphi:\Space R\rightarrow \Space S$
   defined by $\varphi^{\sharp}:S\rightarrow R$,
 then
  $\Space R_{\varphi}$,
    where $R_{\varphi}$ is the image $\varphi^{\sharp}(S)$
     of $S$ in $R$,
  is called the {\it surrogate of $X$ associated to $\varphi$}.
}\end{sdefinition}

\smallskip

As Example~1.4 illustrates,
   % Example [Azumaya point]
 commutative surrogates may be used to manifest/reveal
 the hidden geometry of a noncommutative space.

\smallskip

\begin{sdefinition}
{\bf [push-forward of module].}
{\rm
 Given a morphism $\varphi: X=\Space R\rightarrow Y=\Space S$,
   defined by $\varphi^{\sharp}:S\rightarrow R$,
  and a (left) $R$-module $M$,
 the {\it push-forward} of $M$ from $X$ to $Y$ under $\varphi$,
  in notation $\varphi_{\ast}M$ or $_SM$ when $\varphi$ is understood,
  is defined to be $M$ as a (left) $S$-module via $\varphi^{\sharp}$.
 Since $\Ker\varphi^{\sharp}\cdot M=0$, we say that
  the $S$-module $\varphi_{\ast}M$ on $Y$ is {\it supported}
  on $\varphi(X)\subset Y$.
}\end{sdefinition}

\smallskip

In particular, any $R$-module $M$ on $X=\Space R$
 has a push-forward on any surrogate of $X$.

\bigskip

\begin{flushleft}
{\bf D-branes in an affine noncommutative space
                   \`{a} la Polchinski-Grothendieck Ansatz.}
\end{flushleft}
A {\it D-brane} is geometrically a locus in space-time
 that serves as the boundary condition for open strings.\footnote{This
                              is how one would think of a D-brane
                               to begin with.
                              Later development of string theory enlarges
                               this picture considerably.
                              See [L-Y1: References] to get a glimpse.}
Through this, open strings dictate also the fields and their dynamics
 on D-branes.
In particular, when a collection of D-branes are stacked together,
 the fields on the D-brane that govern the deformation of the brane
 are enhanced to matrix-valued,
 cf.\ Polchinski in [Po: vol.~I, Sec.~8.7].
This open-string-induced phenomenon on D-branes,
 when re-read from Grothendieck's contravariant equivalence
  between the category of geometries and the category of algebras,
 says that D-brane world-volume carries
  an Azumaya-type noncommutative structure.
I.e.
 \begin{itemize}
  \item[$\cdot$]
   {\it Polchinski-Grothendieck Ansatz}$\,$:
   D-brane has a geometry that is generically locally\\
   associated to algebras of the form $M_r(R_0)$,
   where $R_0$ is an ${\Bbb R}$-algebra.
 \end{itemize}
See [L-Y1: Sec.~2.2] for detailed explanations.

For this work, we will be restricting ourselves to
 affine situations in noncommutative algebraic geometry with
 $R_0$ a commutative Noetherian ${\Bbb C}$-algebra.
Thus:

\smallskip

\begin{sdefinition}
{\bf [affine D-brane in affine target].} {\rm
 A {\it D-brane} (or {\it D-brane world-volume})
  in an affine noncommutative space $Y=\Space S$
  is a triple that consists of
  \begin{itemize}
   \item[$\cdot$]
    a ${\Bbb C}$-algebra $R$ that is isomorphic to $M_r(R_0)$
     for an $R_0$,

   \item[$\cdot$]
    a (left) generically simple $R$-module $M$,
     which has rank $r$ as an $R_0$-module,

   \item[$\cdot$]
    a morphism $\varphi:\Space R\rightarrow Y$,
     defined by a ${\Bbb C}$-algebra-homomorphism
     $\varphi^{\sharp}: S\rightarrow R$.
  \end{itemize}
 We will write $\varphi:(\Space R,M)\rightarrow Y$
  for simplicity of notations.
 $\varphi(X)=\Image\varphi$ is called the {\it image-brane} on $Y$.
 $M$ is called the {\it fundamental module} on $\Space R$  and
 the push-forward $\varphi_{\ast}M$ is called
  the {\it Chan-Paton module} on the image-brane $\varphi(X)$.
}\end{sdefinition}

\smallskip

\begin{sdefinition-example}
{\bf [D0-brane as morphism from Azumaya point with fundamental module].}
{\rm
 A {\it D$0$-brane of length $r$}
  on an affine noncommutative space $Y=\Space S$
  is given by a morphism
  $\varphi:(\Space\End(V), V)\rightarrow Y$, where $V\simeq{\Bbb C}^r$.
 In other words, a D$0$-brane on $Y$ is given by
  \begin{itemize}
   \item[$\cdot$]
    a finite-dimensional ${\Bbb C}$-vector space $V$  and
    a ${\Bbb C}$-algebra-homomorphism:
     $\varphi^{\sharp}:S\rightarrow \End(V)$.
  \end{itemize}
  This is precisely a realization of
   a finite-dimensional ${\Bbb C}$-vector space $V$
   as an $S$-module.\footnote{Thus,
                  a D$0$-brane on {\it Space}$\,S$ is precisely
                   an $S$-module
                   that is of finite dimension
                    as a ${\Bbb C}$-vector space.
                  Such a direct realization of a D-brane as a module
                   on a target-space is a special feature
                   for D$0$-branes.
                  For high dimensional D-branes,
                   such modules on the target-space give only
                   a subclass of D-branes
                   that describe solitonic branes in space-time.}
 A {\it morphism}
   from $\varphi_1:(\Space \End(V_1),V_1)\rightarrow Y$
   to   $\varphi_2:(\Space \End(V_2),V_2)\rightarrow Y$
  is a ${\Bbb C}$-vector-space isomorphism
  $h:V_2 \stackrel{\sim}{\rightarrow} V_1$
   such that the following diagram commutes
   $$
    \xymatrix{
       \End(V_1)  &&  S \ar[ll]_{\varphi_1^{\sharp}}
                        \ar[lld]^{\varphi_2^{\sharp}} \\
       \End(V_2) \ar[u]^h   &&&.
    }
   $$
 Here, the $h$-induced isomorphism
  $\End(V_2)\stackrel{\sim}{\rightarrow} \End(V_1)$
  is also denoted by $h$.
 In other words, a morphism between $\varphi_1$ and $\varphi_2$
  is an isomorphism of the corresponding $V_1$ and $V_2$
  as $S$-modules.
}\end{sdefinition-example}

\smallskip

It follows from the above definition/example
       % Definition/Example [D0-brane as morphism from Azumaya point
       %                     with fundamental module]
 that the moduli stack ${\mathfrak M}^{\rm D0}_r(Y)$
  of D$0$-branes of length $r$ on $Y=\Space S$ has an atlas
  given by the {\it representation scheme}
  $\Rep(S,M_r({\Bbb C}))$ that parameterizes
  all ${\Bbb C}$-algebra-homomorphisms $S\rightarrow M_r({\Bbb C})$.
The latter commutative scheme serves also as
 the moduli space of morphisms $\Space M_r({\Bbb C})\rightarrow Y$
 with $M_r({\Bbb C})$ treated as fixed.
{From} [K-R] and [leB1], one expects that
 noncommutative geometric structures/properties of $Y=\Space S$
 are reflected in properties/structures of the discrete family
 of commutative schemes $\Rep(S,M_r({\Bbb C}))\,$, $r\in {\Bbb Z}_{>0}$.
This anticipation from noncommutative algebraic geometry rings
 hand in hand with the stringy philosophy to use D-branes
 as a probe to the nature of space-time!

\bigskip

\section{Deformations of a conifold via an Azumaya probe.}

Using a toric setup for a conifold
  that is meant to match Klebanov-Witten [K-W],
 we discuss how an Azumaya probe ``sees" deformations of the conifold
 in a way that resembles Klebanov-Strassler [K-S].

\bigskip

\begin{flushleft}
{\bf A toric setup for the standard local conifold.}
\end{flushleft}
The standard local conifold
 $\,Y=\Spec ({\Bbb C}[z_1,z_2,z_3,z_4]/(z_1z_2-z_3z_4))\,$
 can be given an affine toric variety description as follows.
Let $N=\oplus_{i=1}^4\,{\Bbb Z}e_i$ be the rank $4$ lattice
 and $\Delta$ be the fan in $N$ that consists of the single
 non-strongly convex polyhedral cone
  $\sigma=\oplus_{i=1}^6 {\Bbb R}_{\ge 0} v_i$
  in $N_{\Bbb R}:=N\otimes_{\Bbb Z}{\Bbb R}$,
  where
  $$
   \begin{array}{l}
    v_1\; =\; e_1\,,\hspace{2em} v_2\; =\; e_2\,,\hspace{2em}
     v_3\; =\; e_3\,,\hspace{2em} v_4\; =\; -e_1+e_2+e_3\,, \\[.6ex]
    v_5\; =\; e_1-e_2-e_3+e_4\,,\hspace{2em}
     v_6\; =\; -v_5\; =\; -e_1+e_2+e_3-e_4\,.
   \end{array}
  $$
Let $M=\Hom(N,{\Bbb Z})$ be the dual lattice of $N$,
 with the dual basis
 $\{e_1^{\ast}\,,\, e_2^{\ast}\,,\, e_3^{\ast}\,,\, e_4^{\ast}\}$.
Then, the dual cone $\sigma^{\vee}$ of $\sigma$ is given by
 $\Span_{{\Bbb R}_{\ge 0}}
  \{\, e_1^{\ast}+e_2^{\ast}\,,\, e_3^{\ast}+e_4^{\ast}\,,\,
       e_1^{\ast}+e_3^{\ast}\,,\, e_2^{\ast}+e_4^{\ast}\, \}\,
  \subset\, M_{\Bbb R}$.
This determines a commutative semigroup
 $$
  S_{\sigma}\;=\; \sigma^{\vee}\cap M\;
   =\; \Span_{{\Bbb Z}_{\ge 0}}
         \{\, e_1^{\ast}+e_2^{\ast}\,,\, e_3^{\ast}+e_4^{\ast}\,,\,
              e_1^{\ast}+e_3^{\ast}\,,\, e_2^{\ast}+e_4^{\ast}\, \}
 $$
 with generators
 $\,e_1^{\ast}+e_2^{\ast}\,,\, e_3^{\ast}+e_4^{\ast}\,,\,
    e_1^{\ast}+e_3^{\ast}\,,\, e_2^{\ast}+e_4^{\ast}\,$.
The corresponding group-algebra
 $$
  {\Bbb C}[S_{\sigma}]\;
   =\; {\Bbb C}[\, \xi_1\xi_2\,,\, \xi_3\xi_4\,,\,
                   \xi_1\xi_3\,,\, \xi_2\xi_4\, ]\;
   \subset\; {\Bbb C}[\xi_1, \xi_2, \xi_3, \xi_4]\,,
 $$
 where $\xi_i=\exp(e_i^{\ast})$, $i=1,\,2,\,3,\,4\,$,
 defines then the conifold
 $$
  Y\; =\; U_{\sigma}\;=\; \Spec ({\Bbb C}[S_{\sigma}])\;
      =\; \Spec({\Bbb C}[z_1,z_2,z_3,z_4]/(z_1z_2-z_3z_4))\,,
 $$
 where
 $$
  z_1\,=\,\xi_1\xi_2\,,\;\; z_2\,=\,\xi_3\xi_4\,,\;\;
  z_3\,=\,\xi_1\xi_3\,,\;\; z_4\,=\,\xi_2\xi_4\,.
 $$
Note that built into this construction is the morphism
 $$
  {\Bbb A}^4_{[\xi_1,\xi_2,\xi_3,\xi_4]}\,
    :=\, \Spec({\Bbb C}[\xi_1,\xi_2,\xi_3,\xi_4])\;
  \longrightarrow\; Y\;
  \hookrightarrow\;
  {\Bbb A}^4_{[z_1,z_2,z_3,z_4]}\,
   :=\, \Spec({\Bbb C}[z_1,z_2,z_3,z_4])\,,
 $$
 where the first morphism is surjective.

\bigskip

\begin{flushleft}
{\bf An Azumaya probe to a noncommutative space and
     its commutative descent.}
\end{flushleft}
Guided by [K-W] and [K-S],
  where $\xi_i$'s here play the role of scalar component
  of chiral superfields involved in ibidem,
consider the noncommutative space
 \begin{eqnarray*}
  \Xi &  :=  & \Space(R_{\Xi}) \\
      &  :=
      & \Space\left(
       \frac{ {\Bbb C}\langle\,\xi_1,\xi_2,\xi_3,\xi_4\, \rangle }
       { ([\xi_1\xi_3, \xi_2\xi_4]\,,\, [\xi_1\xi_3, \xi_1\xi_4]\,,\,
          [\xi_1\xi_3, \xi_2\xi_3]\,,\, [\xi_2\xi_4, \xi_1\xi_4]\,,\,
          [\xi_2\xi_4, \xi_2\xi_3]\,,\, [\xi_1\xi_4, \xi_2\xi_3]) }
              \right)\,,
 \end{eqnarray*}
 where
  ${\Bbb C}\langle \xi_1, \xi_2, \xi_3, \xi_4\rangle$
    is the associative unital ${\Bbb C}$-algebra generated by
    $\{\xi_1,\xi_2,\xi_3,\xi_4\}$,
  $(\,\cdots\,)$ in the denominator is the two-sided ideal
   generated by $\cdots$, and
  $[\,\bullet\,,\,\bullet^{\prime}\,]$ is the commutator.
Here, $\Space(\,\bullet\,)$ is the would-be space associated to
 the ring $\,\bullet\,$.
We do not need its detail as all we need are morphisms between spaces
 which can be contravariantly expressed as ring-homomorphisms.
By construction, the scheme-morphism
 ${\Bbb A}^4_{[\xi_1,\xi_2,\xi_3,\xi_4]}
   \rightarrow {\Bbb A}^4_{[z_1,z_2,z_3,z_4]}$, whose image is $Y$,
 extends to a morphism
 $$
  \pi^{\Xi}\; :\;
    \Xi\; \longrightarrow\; {\Bbb A}^4_{[z_1,z_2,z_3,z_4]}\,,
 $$
 whose image is now the whole ${\Bbb A}^4_{[z_1,z_2,z_3,z_4]}$.
The underlying ring-homomorphism is given by
 $$
 \begin{array}{cccccl}
  \pi^{\Xi,\sharp}  & :
   & {\Bbb C}[z_1,z_2,z_3,z_4]  & \longrightarrow  &  R_{\Xi}\\
  && z_1 & \longmapsto  & \xi_1\xi_3 \\
  && z_2 & \longmapsto  & \xi_2\xi_4 \\
  && z_3 & \longmapsto  & \xi_1\xi_4 \\
  && z_4 & \longmapsto  & \xi_2\xi_3  & .\\
 \end{array}
 $$

Consider a D$0$-brane moving on the conifold $Y$
 via the chiral superfields.
In terms of Polchinski-Grothendieck Ansatz,
 this is realized by the descent of morphisms
 $\widetilde{\varphi}:
   \Space M_1({\Bbb C})=\Spec{\Bbb C}\rightarrow \Xi$
 to $\varphi: \Space M_1({\Bbb C})=\Spec{\Bbb C}\rightarrow Y$
 by the specification of ring-homomorphisms
 $$
  \widetilde{\varphi}^{\sharp}\;:\;
   \xi_1\;\longmapsto\; a_1\,;\hspace{1em}
   \xi_2\;\longmapsto\; a_2\,;\hspace{1em}
   \xi_3\;\longmapsto\; b_1\,;\hspace{1em}
   \xi_4\;\longmapsto\; b_2\,.
 $$
The corresponding
 $$
  \varphi^{\sharp}\;:\;
   z_1\;\longmapsto\; a_1b_1\,;\hspace{1em}
   z_2\;\longmapsto\; a_2b_2\,;\hspace{1em}
   z_3\;\longmapsto\; a_1b_2\,;\hspace{1em}
   z_4\;\longmapsto\; a_2b_1
 $$
 gives a morphism $\varphi:\Spec{\Bbb C}\rightarrow Y$,
 i.e.\ a ${\Bbb C}$-point on the conifold $Y$.

\bigskip

\begin{flushleft}
{\bf Deformations of the conifold via an Azumaya probe:
     descent of noncommutative superficially-infinitesimal deformations.}
\end{flushleft}
We now consider what happens if we add a D$0$-brane to the
conifold point
 of $Y$.
This D$0$-brane together with the D$0$-brane probe
 is the image of a morphism from the Azumaya point
 $\Space M_2({\Bbb C})$ to $Y$.
Thus we should consider morphisms
 $\widetilde{\varphi}:\Space M_2({\Bbb C})\rightarrow \Xi$
 of noncommutative spaces  and
 their descent $\varphi$ on related commutative spaces.

\smallskip

\begin{sdefinition}
{\bf [superficially infinitesimal deformation].} {\rm
 Given finitely-presented associative unital rings,
  $R=\langle\,r_1,\,\ldots\,,r_m\,\rangle/\!\!\sim$ and $S$,  and
  a ring-homomorphism $h:R\rightarrow S$.
 A {\it superficially infinitesimal deformation} of $h$
   {\it with respect to the generators}
   $\{r_1,\,\ldots\,,r_m\}$ {\it of} $R$
  is a ring-homomorphism $h_{\varepsilon}:R\rightarrow S$
   such that
    $h_{\varepsilon}(r_i)=h(r_i)+\varepsilon_i$
     with $\varepsilon_i^2=0$,
    for $i=1,\,\ldots\,,m$.
}
\end{sdefinition}

\smallskip

\begin{sremark}
{\it $[\,$commutative $S$$\,]$.}
{\rm
 Note that when $S$ is commutative,
 a superficially infinitesimal deformation of
  $\, h_{\varepsilon}:R\rightarrow S\,$
  is an infinitesimal deformation of $h$
 in the sense that $h_{\varepsilon}(r)=h(r)+\varepsilon_r$
  with $(\varepsilon_r)^2=0$, for all $r\in R$.
 This is no longer true for general noncommutative $S$.
}
\end{sremark}

\smallskip

To begin, consider the diagram of morphisms of spaces
  $$
  \xymatrix{
   \Space M_2({\Bbb C})
    \ar @{=}[d]\ar[rrrr]^{\widetilde{\varphi}}
    &&&& \Xi=\Space R_{\Xi} \ar[d]^{\pi^{\Xi}} \\
   \Space M_2({\Bbb C})\ar[rrr]^{\hspace{2em}\varphi}
    &&& Y \ar @{^{(}->}[r] & {\Bbb A}^4
  }
  $$
 given by ring-homomorphisms
  $$
  \xymatrix{
   M_2({\Bbb C}) \ar @{=}[d]
    &&&& R_{\Xi} \ar[llll]_{\widetilde{\varphi}^{\sharp}} \\
   M_2({\Bbb C})
    &&& {\Bbb C}[z_1,z_2,z_3,z_4]/(z_1z_2-z_3z_4)
          \ar[lll]_{\varphi^{\sharp}\hspace{3.6em}}
    &   {\Bbb C}[z_1,z_2,z_3,z_4]
          \ar @{->>}[l] \ar[u]_{\pi^{\Xi,\sharp}}
  }
  $$
  with
  $$
  \xymatrix{
   A_1\,;\;\;  A_2\,;\;\;  B_1\,;\;\;  B_2
    &&& \xi_1\,;\;  \xi_2\,;\;  \xi_3\,;\;\;  \xi_4
           \ar@{|->}[lll]_{\widetilde{\varphi}^{\sharp}}
           \\
    &&& \xi_1\xi_3\,;\;\; \xi_2\xi_4\,;\;\;
        \xi_1\xi_4\,;\;\; \xi_2\xi_3 \\
   A_1B_1\,;\;\;  A_2B_2\,;\;\;  A_1B_2\,;\;\;  A_2B_1
    && \overline{z_1}\,;\;\;  \overline{z_2}\,;\;\;
       \overline{z_3}\,;\;\;  \overline{z_4}
        \ar@{|->}[ll]_{\hspace{3.4em}\varphi^{\sharp}}
    &  z_1\,;\;\;  z_2\,;\;\;  z_3\,;\;\;  z_4
        \ar@{|->}[l] \ar@{|->}[u]_{\pi^{\Xi,\sharp}}
  }
  $$
  where
  $$
   A_1\;=\; \left[
             \begin{array}{cc} a_1 & 0 \\ 0 & 0 \end{array}
            \right]\,, \hspace{2em}
   A_2\;=\; \left[
             \begin{array}{cc} a_2 & 0 \\ 0 & 0 \end{array}
            \right]\,, \hspace{2em}
   B_1\;=\; \left[
             \begin{array}{cc} b_1 & 0 \\ 0 & 0 \end{array}
            \right]\,, \hspace{2em}
   B_2\;=\; \left[
             \begin{array}{cc} b_2 & 0 \\ 0 & 0 \end{array}
            \right]\,.
  $$
 The image D-brane $\varphi(\Space M_2({\Bbb C}))$
  is supported on a subscheme $Z$ of $Y$
  associated to the ideal
   $$
    \Ker \varphi\;
    =\;
     \left\{
      \begin{array}{ll}
      (\overline{z_1},\, \overline{z_2},\,
                          \overline{z_3},\, \overline{z_4})\,
          \cap\,
         (\overline{z_1}-a_1b_1,\, \overline{z_2}-a_2b_2,\,
           \overline{z_3}-a_1b_2,\, \overline{z_4}-a_2b_1 )\\[.6ex]
       \hspace{7.2em}\mbox{if the tuple
         $(a_1b_1, a_2b_2, a_1b_2, a_2b_1)\ne(0,0,0,0)$}\,,\\[1.2ex]
      (\overline{z_1},\, \overline{z_2},\,
                         \overline{z_3},\, \overline{z_4} )
       \hspace{1em}\mbox{if the tuple
         $(a_1b_1, a_2b_2, a_1b_2, a_2b_1)=(0,0,0,0)$}\,.
      \end{array}
     \right.
   $$
 The former corresponds to two simple non-coincident D0-branes,
  each with Chan-Paton module ${\Bbb C}$, on the conifold $Y$
  with one of them sitting at the conifold point ${\mathbf 0}$
   and the other sitting at the ${\Bbb C}$-point with
   the coordinate tuple $(a_1b_1,\, a_2b_2,\, a_1b_2,\, a_2b_1)$
 while
 the latter corresponds to coincident D$0$-branes at ${\mathbf 0}$
  with the Chan-Paton module enhanced to ${\Bbb C}^2$
  at ${\mathbf 0}$.
 In both situations, the support $Z$ of the D-brane is reduced.
 This is the transverse-to-the-effective-space-time part
  of the D$3$-brane setting in [K-W] and [K-S].

Consider now a superficially infinitesimal deformation
 of $\widetilde{\varphi}$ given by:
$$
 \begin{array}{c}
 \xymatrix{
  \Space M_2({\Bbb C}) \hspace{1em}
    \ar[rrr]^{\widetilde{\varphi}_{(\delta_1,\delta_2,\eta_1,\eta_2)}}
    &&& \hspace{1em}
        \Xi=\Space R_{\Xi} }\\
 \xymatrix{
  M_2({\Bbb C})\hspace{2em}
    &&& \hspace{3em}
        R_{\Xi}
         \ar[lll]_{\widetilde{\varphi}
                   _{(\delta_1,\delta_2,\eta_1,\eta_2)}^{\sharp}}
         \hspace{1em}  }\\[1.2ex]
 \xymatrix{
  A_1\,;\;\;  A_2\,;\;\;  B_1\,;\;\;  B_2 \hspace{2em}
    && \hspace{2em}
       \xi_1\,;\;  \xi_2\,;\;  \xi_3\,;\;\;  \xi_4
           \ar@{|->}[ll] \hspace{2em} }
 \end{array}
$$
where
 $$
  A_1\;=\; \left[
            \begin{array}{cc} a_1 & \delta_1 \\ 0 & 0 \end{array}
           \right]\,, \hspace{2em}
  A_2\;=\; \left[
            \begin{array}{cc} a_2 & \delta_2 \\ 0 & 0 \end{array}
           \right]\,, \hspace{2em}
  B_1\;=\; \left[
            \begin{array}{cc} b_1 & 0 \\ \eta_1 & 0 \end{array}
           \right]\,, \hspace{2em}
  B_2\;=\; \left[
            \begin{array}{cc} b_2 & 0 \\ \eta_2 & 0 \end{array}
           \right]\,.
 $$
Should $\Space M_2({\Bbb C})$ be a commutative space,
 this would give only an infinitesimal deformation of $\varphi$.
However, $\Space M_2({\Bbb C})$ is not a commutative space
 and, hence, the naive anticipation above could fail.
Indeed, the descent
 $\varphi_{(\delta_1,\delta_2,\eta_1,\eta_2)}$ of
 $\widetilde{\varphi}_{(\delta_1,\delta_2,\eta_1,\eta_2)}$
 is given by
$$
 \begin{array}{l}
 \xymatrix{
  \hspace{3.4em}
  \Space M_2({\Bbb C}) \hspace{4em}
    \ar[rrr]^{\hspace{4em}\varphi_{(\delta_1,\delta_2,\eta_1,\eta_2)}}
    &&& \hspace{3.6em} {\Bbb A}^4 }\\
 \xymatrix{
  \hspace{5em}
  M_2({\Bbb C})\hspace{5em}
    &&& \hspace{1.4em}
        {\Bbb C}[z_1,z_2,z_3,z_4]
         \ar[lll]
          _{\hspace{2em}
            \varphi_{(\delta_1,\delta_2,\eta_1,\eta_2)}^{\sharp}}
         \hspace{1em}  }\\[1.2ex]
 \xymatrix{
  A_1B_1\,;\;\;  A_2B_2\,;\;\;  A_1B_2\,;\;\;  A_2B_1 \hspace{2em}
    && \hspace{2em}
       z_1\,;\;\;  z_2\,;\;\;  z_3\,;\;\;  z_4\,,
           \ar@{|->}[ll] \hspace{2em} }
 \end{array}
$$
i.e.
$$
 \begin{array}{l}
  \mbox{\small
    $\left[\begin{array}{cc}
     a_1b_1+\delta_1\eta_1 & 0 \\ 0 & 0 \end{array}\right]$}\,;\;\;
  \mbox{\small
    $\left[\begin{array}{cc}
     a_2b_2+\delta_2\eta_2 & 0 \\ 0 & 0 \end{array}\right]$}\,;\;\;
  \mbox{\small
    $\left[\begin{array}{cc}
     a_1b_2+\delta_1\eta_2 & 0 \\ 0 & 0 \end{array}\right]$}\,;\;\;
  \mbox{\small
    $\left[\begin{array}{cc}
     a_2b_1+\delta_2\eta_1 & 0 \\ 0 & 0 \end{array}\right]$} \\[3ex]
  \xymatrix{ \hspace{24em}
    && \hspace{2em}
       z_1\,;\;\;  z_2\,;\;\;  z_3\,;\;\;  z_4\,.
           \ar@{|->}[ll] \hspace{2em}
   }
 \end{array}
$$
The image
 $Z:=\varphi_{(\delta_1,\delta_2,\eta_1,\eta_2)}\,
            (\Space M_2({\Bbb C}))$
 of the Azumaya point $\Space M_2({\Bbb C})$
 under $\varphi_{(\delta_1,\delta_2,\eta_1,\eta_2)}$
 remains a $0$-dimensional reduced scheme, consisting of
 either two ${\Bbb C}$-points - with one of them at ${\mathbf 0}$ -
 or ${\mathbf 0}$ alone.
However,
 $$
  z_1z_2-z_3z_4\;=\;
   \left|
    \begin{array}{cc} z_1 & z_3 \\ z_4 & z_2 \end{array}
   \right|\;
  =\;
   \left|
    \begin{array}{cc} a_1 & \delta_1 \\ a_2 & \delta_2 \end{array}
   \right|\,\cdot\,
   \left|
    \begin{array}{cc} b_1 & b_2 \\ \eta_1 & \eta_2 \end{array}
   \right|
 $$
 vanishes if and only if either
  {\footnotesize
   $\left|
    \begin{array}{cc} a_1 & \delta_1 \\ a_2 & \delta_2 \end{array}
   \right|$} or
  {\footnotesize $\left|
    \begin{array}{cc} b_1 & b_2 \\ \eta_1 & \eta_2 \end{array}
   \right|$}
  is $0$.
In other words,
 while the image
  $\varphi_{(\delta_1,\delta_2,\eta_1,\eta_2)}\,(\Space M_2({\Bbb C}))$
 still contains the conifold-point ${\mathbf 0}$ in $Y$,
 as a whole it may longer lie completely even in any infinitesimal
  neighborhood of the conifold $Y$ in ${\Bbb A}^4$.
I.e.:

\smallskip

\begin{slemma}
{\bf [deformation from descent of
      superficially infinitesimal deformation].}
 The descent $\varphi_{(\delta_1,\delta_2,\eta_1,\eta_2)}$
  of a superficially infinitesimal deformation of $\widetilde{\varphi}$
  can truly deform $\varphi$.
 Thus, an appropriate choice of a subspace of the space of morphisms
  $\widetilde{\varphi}_{(\,\bullet\,)}:
                              \Space M_2({\Bbb C})\rightarrow \Xi$
  can descend to give a space of morphisms
  $\varphi_{(\,\bullet\,)}:
      \Space M_2({\Bbb C})\rightarrow {\Bbb A}^4$
  that is parameterized by a deformed conifold $Y^{\prime}$.
\end{slemma}

\smallskip

\noindent
This realizes a deformed conifold
  as a moduli space of morphisms from an Azumaya point  and
 is the reason why the Azumaya probe can see a deformation
  of the conifold $Y$ from the viewpoint of Polchinski-Grothendieck
  Ansatz.
{\sc Figure}~2-1.

\begin{figure}[htbp]
 \epsfig{figure=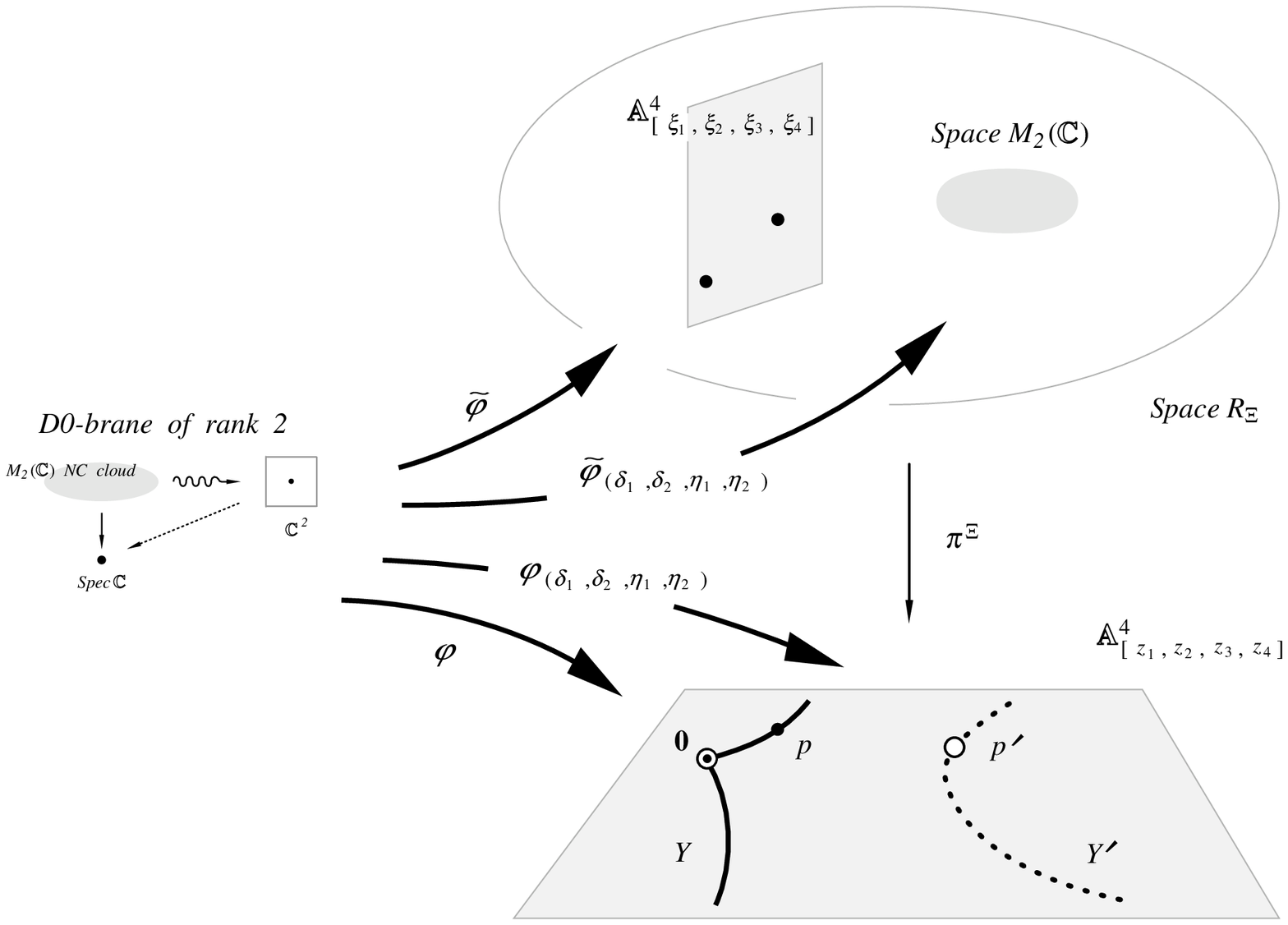,width=16cm}
 \centerline{\parbox{13cm}{\small\baselineskip 12pt
  {\sc Figure}~2-1.
  A generic superficially infinitesimal deformation
   $\widetilde{\varphi}_{(\delta_1,\delta_2,\eta_1,\eta_2)}$
   of $\widetilde{\varphi}$ has
   a noncommutative image $\simeq$ {\it Space}$\,M_2({\Bbb C})$.
  It then descends to ${\Bbb A}^4_{[z_1,z_2,z_3,z_4]}$  and
   becomes a pair of ${\Bbb C}$-points
    on ${\Bbb A}^4_{[z_1,z_2,z_3,z_4]}$.
  One of the points is the conifold singularity
   ${\mathbf 0}=V(z_1,z_2,z_3,z_4)\in Y$  and
  the other is the point
   $p^{\prime}
    =V(\,z_1-a_1b_1-\delta_1\eta_1\,,\, z_2-a_2b_2-\delta_2\eta_2\,,\,
         z_3-a_1b_2-\delta_1\eta_2\,,\, z_4-a_2b_1-\delta_2\eta_1\,)\,$
   {\it off} $\,Y$ (generically).
 Through such deformations, any ${\Bbb C}$-point on
  ${\Bbb A}^4_{[z_1,z_2,z_3,z_4]}$ can be reached.
 Thus, one can realizes a deformation $Y^{\prime}$ of $Y$
   in ${\Bbb A}^4_{[z_1,z_2,z_3,z_4]}$
  by a subvariety in {\it Rep}$\,(R_{\Xi},M_2({\Bbb C}))$.
 This is the Azumaya-geometry origin of
  the phenomenon in Klebanov-Strassler [K-S] that
   a trapped D-brane sitting on the conifold singularity
   may give rise to a deformation of the moduli space of
   SQFT on the D$3$-brane probe, turning a conifold to
   a deformed conifold.
 Our D$0$-brane here corresponds to the internal part of
  the effective-space-time-filling D$3$-brane world-volume of [K-S].
 }}
\end{figure}

\smallskip

\begin{sremark}
{\it $[\,$generalization$\,]$.}
{\rm
 This phenomenon can be generalized beyond a conifold.
 In particular,
  recall that an $A_n$-singularity
   on a complex surface is also a toric singularity.
 Similar mechanism/discussion can be applied to deform
  a transverse $A_n$-singularity via morphisms from
  an Azumaya probe.
 }
\end{sremark}

\smallskip

\bigskip

\begin{flushleft}
{\bf  Deformations of the conifold via an Azumaya probe: details.}
\end{flushleft}
We now give an explicit construction that realizes Lemma~2.3.
        % Lemma [deformation from descent of superficially
        %        infinitesimal deformation       ]
For convenience\footnote{If
                    {\it Space}$\,M_2({\Bbb C})$ is not fixed,
                     then one studies Artin stacks that
                     parameterizes morphisms in question
                     from {\it Space}$\,M_2({\Bbb C})$ to
                     {\it Space}$\,R_{\Xi}$, the conifold $Y$,
                     and ${\Bbb A}^4_{[z_1,z_2,z_3,z_4]}$
                     respectively.
                    The discussion given here is then
                     on an atlas of the stack in question.},
 we will take $\Space M_2({\Bbb C})$ as fixed,
 and is equipped with
 the defining fundamental (left) $M_2({\Bbb C})$-module ${\Bbb C}^2$.
Then, the space $\Mor^a(\Space M_2({\Bbb C}),\Xi)$ of admissible morphisms
 of the form $\widetilde{\varphi}_{(\,\bullet\,)}$ in the previous theme
 is naturally realized as a subscheme $\Rep^a(R_{\Xi},M_2({\Bbb C}))$
 of the representation scheme $\Rep(R_{\Xi},M_2({\Bbb C}))$
 that parameterizes elements in
 $\Mor_{{\Bbb C}\mbox{\scriptsize-}\Algscriptsize}(R_{\Xi}, M_2({\Bbb C}))$.
{From} the previous discussion,
 \begin{eqnarray*}
  \lefteqn{\Rep^a(R_{\Xi},M_2({\Bbb C}))\;
    =\; \Spec {\Bbb C}[a_1,a_2,\delta_1,\delta_2, b_1,b_2,\eta_1,\eta_2]
        }\\[.6ex]
  && =:\;
     {\Bbb A}^8_{[a_1,a_2,\delta_1,\delta_2,b_1,b_2,\eta_1,\eta_2]}\;
     =\; {\Bbb A}^4_{[a_1,a_2,\delta_1,\delta_2]}
      \times_{\Bbb C} {\Bbb A}^4_{[b_1,b_2,\eta_1,\eta_2]}\,.
 \end{eqnarray*}
Consider also the space $\Mor^a(\Space M_2({\Bbb C}),{\Bbb A}^4)$
 of morphisms from Azumaya point to ${\Bbb A}^4_{[z_1,z_2,z_3,z_4]}$
 with the associated ${\Bbb C}$-algebra-homomorphism of the form
 $$
  z_1\;\longmapsto\;
   \left[\begin{array}{cc} c_1 & 0 \\ 0 & 0 \end{array}\right]\,,\;\;
  z_2\;\longmapsto\;
   \left[\begin{array}{cc} c_2 & 0 \\ 0 & 0 \end{array}\right]\,,\;\;
  z_3\;\longmapsto\;
   \left[\begin{array}{cc} c_3 & 0 \\ 0 & 0 \end{array}\right]\,,\;\;
  z_4\;\longmapsto\;
   \left[\begin{array}{cc} c_4 & 0 \\ 0 & 0 \end{array}\right]\,.
 $$
Denote the associated representation scheme by
 $$
  \Rep^a({\Bbb C}[z_1,z_2,z_3,z_4], M_2({\Bbb C}))\,,
  \hspace{1em}\mbox{which is}\hspace{1em}
  \Spec{\Bbb C}[c_1,c_2,c_3,c_4]\;=:\; {\Bbb A}^4_{[c_1,c_2,c_3,c_4]}\,.
 $$

The ${\Bbb C}$-algebra homomorphism
 $\pi^{\Xi,\sharp}:{\Bbb C}[z_1,z_2,z_3,z_4]\rightarrow R_{\Xi}$
 induces a morphism of representation schemes
 $$
  \pi_{\scriptsizeRep}\;:\;
   {\Bbb A}^8_{[a_1,a_2,\delta_1,\delta_2,b_1,b_2,\eta_1,\eta_2]}\;
   \longrightarrow\; {\Bbb A}^4_{[c_1,c_2,c_3,c_4]}
 $$
 with $\pi_{\scriptsizeRep}^{\sharp}$ given in a matrix form by
 $$
  \pi_{\scriptsizeRep}^{\sharp}\;:\;
  \left[\begin{array}{cc} c_1 & c_3 \\ c_4 & c_2 \end{array}\right]\;
  \longmapsto\;
      \left[
       \begin{array}{cc} a_1 & \delta_1 \\ a_2 & \delta_2 \end{array}
      \right]\,\cdot\,
      \left[
       \begin{array}{cc} b_1 & b_2 \\ \eta_1 & \eta_2 \end{array}
      \right]\,.
 $$

\smallskip

\begin{slemma}
{\bf [enough superficially infinitesimally deformed morphisms].}
  $$
   \pi_{\scriptsizeRep}\;:\;
    {\Bbb A}^8_{[a_1,a_2,\delta_1,\delta_2,b_1,b_2,\eta_1,\eta_2]}\;
    \longrightarrow\; {\Bbb A}^4_{[c_1,c_2,c_3,c_4]}
  $$
  is surjective.
\end{slemma}

\smallskip

There are three homeomorphism classes of fibers of
 $\pi_{\scriptsizeRep}$ over a closed point of
 ${\Bbb A}^4_{[c_1,c_2,c_3,c_4]}$, depending on the rank of
 {\footnotesize
  $\left[\begin{array}{cc} c_1 & c_3 \\ c_4 & c_2 \end{array}\right]$}.

\smallskip

\begin{slemma}
{\bf [topological type of fibers of $\pi_{\scriptsizeRep}$].}
 Let $C^3_{[c_1,c_2,c_3,c_4]}$ be the subvariety of
  ${\Bbb A}^4_{[c_1,c_2,c_3,c_4]}$ associated to the ideal
  $(c_1c_2-c_3c_4)$.
 Similarly, for $C^3_{[a_1,a_2,\delta_1,\delta_2]}$ and
  $C^3_{[b_1,b_2,\eta_1,\eta_2]}$.
 Then:
  \begin{itemize}
   \item[(0)]
    Over ${\mathbf 0}$,
    the fiber is given by
    ${\Bbb A}^4_{[a_1, a_2,\delta_1, \delta_2]}
     \cup {\Bbb A}^4_{[b_1, b_1, \eta_1, \eta_2]}
     \cup \Pi^5$,
    where $\Pi^5$ is a $5$-dimensional irreducible affine scheme
     meeting
     ${\Bbb A}^4_{[a_1, a_2, \delta_1, \delta_2]}
       \cup {\Bbb A}^4_{[b_1, b_2, \eta_1, \eta_2]}$
     along
     $C^3_{[a_1, a_2\, \delta_1, \delta_2]}
      \cup  C^3_{[b_1, b_2\, \eta_1, \eta_2]}$.

   \item[(1)]
    Over a closed point of $C^3_{[c_1,c_2,c_3,c_4]}-\{\mathbf{0}\}$,
    the fiber is the union $\Pi^4_1\cup \Pi^4_2$ of
    two irreducible $4$-dimensional affine scheme
     meeting at a deformed conifold.

   \item[(2)]
    Over a closed point of
     ${\Bbb A}^4_{[c_1,c_2,c_3,c_4]}-C^3_{[c_1,c_2,c_3,c_4]}$,
    the fiber is isomorphic to
     ${\Bbb A}^4_{[a_1, a_2\, \delta_1, \delta_2]}
      - C^3_{[a_1, a_2\, \delta_1, \delta_2]}
     \simeq
     {\Bbb A}^4_{[b_1, b_2\, \eta_1, \eta_2]}
      - C^3_{[b_1, b_2\, \eta_1, \eta_2]}$.
   \end{itemize}
\end{slemma}

\bigskip

\noindent
The lemma follows from a straightforward computation.\footnote{It
                         is very instructive to think of the fibration
                         $\pi_{Rep}:
                         {\Bbb A}^8 _{[a_1,a_2,\delta_1,\delta_2,
                                       b_1,b_2,\eta_1,\eta_2]}
                          \rightarrow {\Bbb A}^4_{[c_1,c_2,c_3,c_4]}$
                         as defining a one-matrix-parameter family
                         of ``matrix nodal curves"
                         in the sense of noncommutative geometry.}
Note that the fundamental group as an analytic space is given by
 $$
 \begin{array}{ccccc}
  \pi_1({\Bbb A}^4_{[c_1, c_2, c_3, c_4]}
                             -C^3_{[c_1, c_2, c_3, c_4]})
   & \simeq
   & \pi_1( {\Bbb A}^4_{[a_1, a_2, \delta_1, \delta_2]}
                - C^3_{[a_1, a_2, \delta_1, \delta_2]})  \\[.6ex]
  &\simeq
   & \pi_1({\Bbb A}^4_{[b_1, b_2, \eta_1, \eta_2]}
                         - C^3_{[b_1, b_2, \eta_1, \eta_2]})
   & \simeq  & {\Bbb Z}
 \end{array}
 $$
 and that 
the smooth bundle-morphism
 $$
  \pi_{\scriptsizeRep}\; :\;
   {\Bbb A}^8_{[a_1,a_2,\delta_1,\delta_2,b_1,b_2,\eta_1,\eta_2]}
    - \pi_{\scriptsizeRep}^{-1}(C^3_{[c_1, c_2, c_3, c_4]})\;
   \longrightarrow\;
   {\Bbb A}^4_{[c_1, c_2, c_3, c_4]} - C^3_{[c_1, c_2, c_3, c_4]}
 $$
 exhibits a monodromy behavior which resembles that of a Dehn twist.

The map
 $\pi_{\scriptsizeRep}:
  {\Bbb A}^8_{[a_1,a_2,\delta_1,\delta_2,b_1,b_2,\eta_1,\eta_2]}
  \rightarrow
  {\Bbb A}^4_{[c_1, c_2, c_3, c_4]}$
admits sections, i.e.\
 morphism
  $s: {\Bbb A}^4_{[c_1, c_2, c_3, c_4]} \rightarrow
      {\Bbb A}^8_{[a_1,a_2,\delta_1,\delta_2,b_1,b_2,\eta_1,\eta_2]}$
  such that $\pi_{\scriptsizeRep}\circ s =$
   the identity map on ${\Bbb A}^4_{[c_1, c_2, c_3, c_4]}$.

\bigskip

\begin{sexample}
{\bf [section of $\pi_{\scriptsizeRep}$].}
{\rm
 Let $t\in \GL_2({\Bbb C})$, then a simple family of sections
  of $\pi_{\scriptsizeRep}$
 $$
  s_{t}\; :\;
    {\Bbb A}^4_{[c_1, c_2, c_3, c_4]}\; \longrightarrow\;
       {\Bbb A}^8_{[a_1,a_2,\delta_1,\delta_2,b_1,b_2,\eta_1,\eta_2]}
 $$
 is given compactly in a matrix expression by
 (with $t$ also in its defining $2\times 2$-matrix form)
 $$
  s_t^{\sharp}\;:\;
   \left(
    \left[ \begin{array}{cc} a_1 & \delta_1  \\ a_2 & \delta_2
     \end{array} \right]\, ,\,
    \left[ \begin{array}{cc} b_1 & b_2  \\ \eta_1 & \eta_2
           \end{array} \right]
   \right)\;
   \longmapsto\;
   \left(
    \left[ \begin{array}{cc} c_1 & c_3 \\ c_4 & c_2 \end{array}
     \right]\,
    \cdot\, t^{-1}\,,\, t\,
   \right)\,.
 $$}
\end{sexample}

\bigskip

\noindent
Through any section
 $s: {\Bbb A}^4_{[c_1, c_2, c_3, c_4]} \rightarrow
     {\Bbb A}^8_{[a_1,a_2,\delta_1,\delta_2,b_1,b_2,\eta_1,\eta_2]}$,
 one can realize $Y^{\prime}\amalg \{\mathbf 0\}$,
  where
   $Y^{\prime}$ is a deformation of the conifold $Y$
    in ${\Bbb A}^4={\Bbb A}^4_{[z_1, z_2, z_3, z_4]}$ and
   ${\mathbf 0}$ is the singular point on $Y$,
 as the descent of a family of superficially infinitesimal deformations
  of morphisms from Azumaya point to the noncommutative space $\Xi$.
In string theory words,
 \begin{itemize}
  \item[$\cdot$] {\it
   deformations of a conifold via a D-brane probe
    are realized by turning on D-branes at the singularity
    appropriately;
   the conifold is deformed and becomes smooth
    while leaving the trapped D-branes at the singularity behind.}
 \end{itemize}
Cf.~{\sc Figure}~2-1.

\bigskip

\section{Resolutions of a conifold via an Azumaya probe.}

In this section, we consider resolutions of
 the conifold $Y=\Spec({\Bbb C}[z_1, z_2, z_3, z_4]/(z_1z_2-z_3z_4))$
 from the viewpoint of an Azumaya probe.
Recall the following diagram of resolutions of $Y$ from blow-ups of $Y$:
 $$
  \xymatrix{
   & \widetilde{Y}
      \ar[ld]_{f_+}
      \ar[dd]^{\pi}
      \ar[rd]^{f_-} & \\
   Y_+\ar[rd]_{\pi_+} && Y_-\ar[ld]^{\pi_-}\\
   & Y &&,
   }
 $$
where
 \begin{itemize}
  \item[$\cdot$]
   $\pi:\widetilde{Y}
      =\Bl_{V(I)}Y=\Proj(\oplus_{i=0}^{\infty}\,I^i)
       \rightarrow Y$
     with $I=(z_1, z_2, z_3, z_4)$,

  \item[$\cdot$]
   $\pi_+:Y_+=\Bl_{V(I_+)}Y=\Proj(\oplus_{i=0}^{\infty}\,I_+^i)
     \rightarrow Y$ with $I_+=(z_1, z_3)$,  and

  \item[$\cdot$]
   $\pi_-: Y_-=\Bl_{V(I_-)}Y=\Proj(\oplus_{i=0}^{\infty}\,I_-^i)
    \rightarrow Y$ with $I_-=(z_1, z_4)$
 \end{itemize}
 are blow-ups of $Y$ along the specified subschemes $V(\,\bullet\,)$
 associated respectively to the ideals $I$, $I_+$, and $I_-$ of
 ${\Bbb C}[z_1, z_2, z_3, z_4]/(z_1z_2-z_3z_4)$ as given.
Here, we set
 $I_{(\pm)}^0={\Bbb C}[z_1, z_2, z_3, z_4]/(z_1z_2-z_3z_4)$.
Let ${\mathbf 0}=V(z_1, z_2, z_3, z_4)$ be the singular point of $Y$.
Then the exceptional locus in each case is given respectively by
 $\pi^{-1}({\mathbf 0})\simeq {\Bbb P}^1\times {\Bbb P}^1$,
 $\pi_+^{-1}({\mathbf 0})\simeq {\Bbb P}^1$, and
 $\pi_-^{-1}({\mathbf 0})\simeq {\Bbb P}^1$;
$Y_+$ and $Y_-$ as schemes/$Y$ are related by a flop;  and
the restriction of birational morphisms
 $f_{\pm}:\widetilde{Y}\rightarrow Y_{\pm}$ to $\pi^{-1}({\mathbf 0})$
 corresponds to the projections of ${\Bbb P}^1\times {\Bbb P}^1$
 to each of its two factors.

\bigskip

\begin{flushleft}
{\bf D-brane probe resolutions of a conifold via the Azumaya structure.}
\end{flushleft}
An atlas for the stack of morphisms from $\Space M_2({\Bbb C})$ to $Y$
 is given by the representation scheme
 $\Rep({\Bbb C}[z_1, z_2, z_3, z_4]/(z_1z_2-z_3z_4), M_2({\Bbb C}))$
 with the $\PGL_2({\Bbb C})$-action
 induced from the $\GL_2({\Bbb C})$-action
 on the fundamental module ${\Bbb C}^2$.
For convenience, we will also call this a $\GL_2({\Bbb C})$-action
 on $\Rep({\Bbb C}[z_1, z_2, z_3, z_4]/(z_1z_2-z_3z_4), M_2({\Bbb C}))$.
Let
 $$
  W\;=\;
  \Rep^{\singletonscriptsize}
  ({\Bbb C}[z_1, z_2, z_3, z_4]/(z_1z_2-z_3z_4), M_2({\Bbb C}))
 $$
 be the subscheme of
 $\Rep({\Bbb C}[z_1, z_2, z_3, z_4]/(z_1z_2-z_3z_4), M_2({\Bbb C}))$
 that parameterizes D$0$-branes
 $\varphi:(\Spec{\Bbb C}, M_2({\Bbb C}), {\Bbb C}^2)\rightarrow Y$
 with $(\Image\varphi)_{\redscriptsize}$
 a single ${\Bbb C}$-point on $Y$.
Explicitly,
 $W$ is the image scheme of
 $$
  \GL_2({\Bbb C})\times W_{ut}\;
   \longrightarrow\;
   \Rep({\Bbb C}[z_1, z_2, z_3, z_4]/(z_1z_2-z_3z_4), M_2({\Bbb C}))
 $$
 where
 $$
  \begin{array}{ccl}
   W_{ut}  & =
    & \left\{
        \rho : {\Bbb C}[z_1, z_2, z_3, z_4]/(z_1z_2-z_3z_4)
                                  \rightarrow M_2({\Bbb C})\;
          \left|\rule{0em}{1.4em}\right.\;
    \mbox{$\rho(z_i)$ is of the form
      $\left[
       \begin{array}{cc} a_i & \varepsilon_i \\ 0  & a_i \end{array}
       \right]$}\,
   \right\}  \\[2.4ex]
  &  \subset
   & \Rep({\Bbb C}[z_1, z_2, z_3, z_4]/(z_1z_2-z_3z_4), M_2({\Bbb C}))
 \end{array}
 $$
 and the morphism $\longrightarrow$ is from the restriction of
  the $\GL_2({\Bbb C})$-group on
  $\Rep({\Bbb C}[z_1, z_2, z_3, z_4]/(z_1z_2-z_3z_4), M_2({\Bbb C}))$.
Using this notation, as a scheme,
 $$
  \begin{array}{ccll}
   W_{ut} & =
     & \Spec ({\Bbb C}[a_1, a_2, a_3, a_4,  \varepsilon_1,
                    \varepsilon_2, \varepsilon_3, \varepsilon_4]
               /(\,a_1a_2-a_3a_4\,,\,
                   a_2\varepsilon_1 + a_1\varepsilon_2
                   - a_4\varepsilon_3 - a_3\varepsilon_4\,) ) \\[1.2ex]
   & \subset
     & \Spec(
         {\Bbb C}[a_1, a_2, a_3, a_4,  \varepsilon_1,
           \varepsilon_2, \varepsilon_3, \varepsilon_4])\;\;\;
       =:\;\;\; {\Bbb A}^{8}_{[a_1, a_2, a_3, a_4,  \varepsilon_1,
                      \varepsilon_2, \varepsilon_3, \varepsilon_4]} &.
  \end{array}
 $$

Imposing the trivial $\GL_2({\Bbb C})$-action on $Y$,
then by construction,
 there is a natural $\GL_2({\Bbb C})$-equivariant morphism
  $$
   \pi^W\; :\; W\; \longrightarrow\;  Y
  $$
  defined by
  $\pi^{W,\sharp}(z_i) = \frac{1}{2}\Tr\rho(z_i) =a_i$
  in the above notation.
This is the morphism that sends
 a $\varphi:(\Spec{\Bbb C}, M_2({\Bbb C}), {\Bbb C}^2)\rightarrow Y$
  under study
 to $(\Image\varphi)_{\redscriptsize}\in Y$.

\bigskip

\begin{slemma}
{\bf [Azumaya probe to conifold singularity].}
 There exists $\GL_2({\Bbb C})$-invariant subschemes
   $\widetilde{Y}^{\prime}$, $Y_+^{\prime}$, and $Y_-^{\prime}$
   of $W$
  such that
   their geometric quotient
    $\widetilde{Y}^{\prime}/\GL_2({\Bbb C})$,
    $Y_+^{\prime}/\GL_2({\Bbb C})$, $Y_-^{\prime}/\GL_2({\Bbb C})$
    under the $\GL_2({\Bbb C})$-action exist
     and are isomorphic to
    $\widetilde{Y}$, $Y_+$, and $Y_-$ respectively.
 Furthermore,
  under these isomorphisms,
   the restriction of $\pi^W: W\rightarrow Y$
   to $\widetilde{Y}^{\prime}$, $Y_+^{\prime}$, and $Y_-^{\prime}$
    descends to morphisms from the quotient spaces
    $\widetilde{Y}^{\prime}/\GL_2({\Bbb C})$,
    $Y_+^{\prime}/\GL_2({\Bbb C})$, $Y_-^{\prime}/\GL_2({\Bbb C})$
    to $Y$
    that realize the resolution diagram
  $$
   \xymatrix{
    & \widetilde{Y}
       \ar[ld]_{f_+}
       \ar[dd]^{\pi}
       \ar[rd]^{f_-} & \\
    Y_+\ar[rd]_{\pi_+} && Y_-\ar[ld]^{\pi_-}\\
    & Y
    }
  $$
  of $Y$ at the beginning of this section.
\end{slemma}

\bigskip

It is in the sense of the above lemma we say that
 \begin{itemize}
  \item[$\cdot$]
   {\it an Azumaya point of rank $\ge 2$ and hence
     a D-brane probe of multiplicity $\ge 2$ can ``see"
     all the three different resolutions of the conifold singularity.}
 \end{itemize}
It should also be noted that Lemma~3.1
             % Lemma [Azumaya probe to conifold singularity]
 is a special case of
 a more general statement that reflects the fact that
 the stack of morphisms from Azumaya points to
 a (general, possibly singular, Noetherian) scheme $Y$
 is a generalization of the notion of jet-schemes of $Y$.
Cf.~[L-Y2: Figure~0-1, caption].

\bigskip

\begin{flushleft}
{\bf An explicit construction of
     $\widetilde{Y}^{\prime}$, $Y_+^{\prime}$, and $Y_-^{\prime}$.}
\end{flushleft}
An explicit construction of
 $\widetilde{Y}^{\prime}$, $Y_+^{\prime}$, and $Y_-^{\prime}$,
 and hence the proof of Lemma~3.1,
                      % Lemma [Azumaya probe to conifold singularity]
 follows from a lifting-to-$W$ of an affine atlas of
 $\Proj(\oplus_{i=0}^{\infty}\,I_{(\pm)}^i)$.

To construct $\widetilde{Y}^{\prime}$,
 recall that $I=(z_1, z_2, z_3, z_4)$.
An affine atlas of $\widetilde{Y}$ is given by the collection
 $$
  U^{(z_i)} =
   \Spec( (\oplus_{j=0}^{\infty}\,I^j)[z_i^{-1}]_0 )\;
   \simeq \left\{
   \begin{array}{ll}
    \Spec( {\Bbb C}[z_1, u_2, u_3, u_4]/(u_2-u_3u_4))
     \simeq {\Bbb A}^3_{[z_1, u_3, u_4]} & \mbox{for $i=1$}\,; \\[1.2ex]
    \Spec( {\Bbb C}[u_1, z_2, u_3, u_4]/(u_1-u_3u_4))
     \simeq {\Bbb A}^3_{[z_2, u_3, u_4]} & \mbox{for $i=2$}\,; \\[1.2ex]
    \Spec( {\Bbb C}[u_1, u_2, z_3, u_4]/(u_1u_2-u_4))
     \simeq {\Bbb A}^3_{[u_1, u_2, z_3]} & \mbox{for $i=3$}\,; \\[1.2ex]
    \Spec( {\Bbb C}[u_1, u_2, u_3, z_4]/(u_1u_2-u_3))
     \simeq {\Bbb A}^3_{[u_1, u_2, z_4]} & \mbox{for $i=4$}\,.
   \end{array} \right.
 $$
Here,
 $z_i\in I$ has grade $1$  and
 $(\oplus_{j=0}^{\infty}\,I^j)[z_i^{-1}]_0$ is the grade-$0$ component
  of the graded algebra $(\oplus_{j=0}^{\infty}\,I^j)[z_i^{-1}]$.
Each $U^{(z_i)}$ is equipped with a built-in morphism
  $\pi^{(i)}:U^{(z_i)}\rightarrow Y$
 in such a way that, when all four are put together,
  they glue to give the resolution $\pi:\widetilde{Y}\rightarrow Y$.

Consider the lifting
  $\{\pi^{(i)\,\prime}: U^{(z_i)} \rightarrow W\, |\,
                                    i\,=\,1\,, 2\,, 3\,, 4\,\}$
  of the atlas
  $\{\, \pi^{(i)}:U^{(z_i)}\rightarrow Y\,|\, i\,=\,1\,, 2\,, 3\,, 4\,\}$
  of $\widetilde{Y}$
 that is given by the lifting
  $\{\pi^{(i)\,\prime}:U^{(z_i)} \rightarrow W_{ut}\subset W\, |\,
                                    i\,=\,1\,, 2\,, 3\,, 4\,\}$
  defined by

  {\small
  $$
   \begin{array}{ccccc}
    \pi^{(1)\,\prime,\sharp}   & :
      & a_1\,,\, a_2\,,\, a_3\,,\, a_4\,,\,
        \varepsilon_1\,,\, \varepsilon_2\,,\,
        \varepsilon_3\,,\, \varepsilon_4
      & \longmapsto
      &  z_1\,,\, z_1u_2\,,\, z_1u_3\,,\, z_1u_4\,,\,
           1\,,\,    u_2\,,\,    u_3\,,\,    u_4
       \hspace{1em}\mbox{respectively}\,,\\[.6ex]
    \pi^{(2)\,\prime,\sharp}   & :
      & a_1\,,\, a_2\,,\, a_3\,,\, a_4\,,\,
        \varepsilon_1\,,\, \varepsilon_2\,,\,
        \varepsilon_3\,,\, \varepsilon_4
      & \longmapsto
      & z_2u_1\,,\,  z_2\,,\, z_2u_3\,,\,  z_2u_4\,,\,
           u_1\,,\,    1\,,\,    u_3\,,\,     u_4
        \hspace{1em}\mbox{respectively}\,,\\[.6ex]
    \pi^{(3)\,\prime,\sharp}   & :
      & a_1\,,\, a_2\,,\, a_3\,,\, a_4\,,\,
        \varepsilon_1\,,\, \varepsilon_2\,,\,
        \varepsilon_3\,,\, \varepsilon_4
      & \longmapsto
      & z_3u_1\,,\,  z_3u_2\,,\, z_3\,,\, z_3u_4\,,\,
           u_1\,,\,     u_2\,,\,   1\,,\,    u_4
        \hspace{1em}\mbox{respectively}\,,\\[.6ex]
    \pi^{(4)\,\prime,\sharp}   & :
      & a_1\,,\, a_2\,,\, a_3\,,\, a_4\,,\,
        \varepsilon_1\,,\, \varepsilon_2\,,\,
        \varepsilon_3\,,\, \varepsilon_4
      & \longmapsto
      & z_4u_1\,,\, z_4u_2\,,\, z_4u_3\,,\, z_4\,,\,
           u_1\,,\,    u_2\,,\,    u_3\,,\,   1
        \hspace{1em}\mbox{respectively}\,.
   \end{array}
  $$
  }

\noindent
$\pi^{(i)\,\prime}$, $i=1,\,2,\,3,\,4$, are now embeddings into $W$
 with the property that
  for any geometric point
   $p\in U^{(z_i)}\times_{\widetilde{Y}}U^{(z_j)}$,
   $\pi^{(i)\,\prime}(p)$ and $\pi^{(j)\,\prime}(p)$
   lies in the same $\GL_2({\Bbb C})$-orbit in $W$.
In other words, up to the pointwise $\GL_2({\Bbb C})$-action,
 they are gluable.
Let $\widetilde{Y}^{\prime}$ be the image scheme of the morphism
 $$
  \GL_2({\Bbb C}) \times
   (U^{(z_1)} \amalg U^{(z_2)} \amalg U^{(z_3)} \amalg U^{(z_4)})\;
  \longrightarrow\; W
 $$
 via $\pi^{(1)\,\prime} \amalg \pi^{(2)\,\prime}
      \amalg \pi^{(3)\,\prime} \amalg \pi^{(4)\,\prime}$
  and the $\GL_2({\Bbb C})$-action on $W$.
Then it follows that
 the geometric quotient $\widetilde{Y}^{\prime}/\GL_2({\Bbb C})$ exists
  and is equipped with a built-in isomorphism
  $\widetilde{Y}^{\prime}/\GL_2({\Bbb C})
                      \stackrel{\sim}{\rightarrow} \widetilde{Y}$,
  as schemes over $Y$,
  through the defining embeddings
   $U^{(z_i)}\hookrightarrow \widetilde{Y}$, $i=1,\, 2,\, 3,\, 4\,$.

\bigskip

For $Y_+^{\prime}$, recall that $I_+=(z_1, z_3)$.
An affine atlas of $Y_+$ is given by the collection
 $$
  U_+^{(z_i)} =
   \Spec( (\oplus_{j=0}^{\infty}\,I^j)[z_i^{-1}]_0 )\;
   \simeq \left\{
   \begin{array}{ll}
    \Spec( {\Bbb C}[z_1, z_2, u_3, z_4]/(z_2-z_4u_3))
     \simeq {\Bbb A}^3_{[z_1, u_3, z_4]} & \mbox{for $i=1$}\,; \\[1.2ex]
    \Spec( {\Bbb C}[u_1, z_2, z_3, z_4]/(z_2u_1-z_4))
     \simeq {\Bbb A}^3_{[u_1, z_2, z_3]} & \mbox{for $i=3$}\,.
   \end{array} \right.
 $$
Each $U_+^{(z_i)}$ is equipped with a built-in morphism
  $\pi_+^{(i)}:U_+^{(z_i)}\rightarrow Y$
 in such a way that, when both are put together,
  they glue to give the resolution $\pi_+:Y_+\rightarrow Y$.

Consider the lifting
  $\{\pi_+^{(i)\,\prime}: U_+^{(z_i)} \rightarrow W\, |\,
                                      i\,=\,1\,, 3\,\}$
  of the atlas
  $\{\, \pi_+^{(i)}:U_+^{(z_i)}\rightarrow Y\,|\, i\,=\,1\,, 3\,\}$
  of $Y_+$
 that is given by the lifting
  $\{\pi_+^{(i)\,\prime}:U_+^{(z_i)} \rightarrow W_{ut}\subset W\, |\,
                                     i\,=\,1\,, 3\,\}$
  defined by

  {\small
  $$
   \begin{array}{ccccc}
    \pi_+^{(1)\,\prime,\sharp}   & :
      & a_1\,,\, a_2\,,\, a_3\,,\, a_4\,,\,
        \varepsilon_1\,,\, \varepsilon_2\,,\,
        \varepsilon_3\,,\, \varepsilon_4
      & \longmapsto
      & z_1\,,\, z_4u_3\,,\, z_1u_3\,,\, z_4\,,\,
          1\,,\,      0\,,\,    u_3\,,\,   0
       \hspace{1em}\mbox{respectively}\,,\\[.6ex]
    \pi_+^{(3)\,\prime,\sharp}   & :
      & a_1\,,\, a_2\,,\, a_3\,,\, a_4\,,\,
        \varepsilon_1\,,\, \varepsilon_2\,,\,
        \varepsilon_3\,,\, \varepsilon_4
      & \longmapsto
      & z_3u_1\,,\, z_2\,,\, z_3\,,\, z_2u_1\,,\,
           u_1\,,\,   0\,,\,   1\,,\,      0
        \hspace{1em}\mbox{respectively}\,.
   \end{array}
  $$
  }

\noindent
The pair,  $\pi_+^{(1)\,\prime}$ and $\pi_+^{(3)\,\prime}$,
 are now embeddings into $W$ that, as in the case of $\widetilde{Y}$,
 are gluable up to the pointwise $\GL_2({\Bbb C})$-action.
Same construction as in the case of $\widetilde{Y}$ gives then
 a $\GL_2({\Bbb C})$-invariant subscheme $Y_+^{\prime}$ of $W$
 whose geometric quotient $Y_+^{\prime}/\GL_2({\Bbb C})$ exists
  and is equipped with a built-in isomorphism
  $Y_+^{\prime}/\GL_2({\Bbb C})\stackrel{\sim}{\rightarrow} Y_+$
 as schemes over $Y$.

\bigskip

For $Y_-^{\prime}$, recall that $I_-=(z_1, z_4)$.
The construction is identical to that in the case of $Y_+$
 after relabelling.
An affine atlas of $Y_-$ is given by the collection
 $$
  U_-^{(z_i)} =
   \Spec( (\oplus_{j=0}^{\infty}\,I^j)[z_i^{-1}]_0 )\;
   \simeq \left\{
   \begin{array}{ll}
    \Spec( {\Bbb C}[z_1, z_2, z_3, u_4]/(z_2-z_3u_4))
     \simeq {\Bbb A}^3_{[z_1, z_3, u_4]} & \mbox{for $i=1$}\,; \\[1.2ex]
    \Spec( {\Bbb C}[u_1, z_2, z_3, z_4]/(z_2u_1-z_3))
     \simeq {\Bbb A}^3_{[u_1, z_2, z_4]} & \mbox{for $i=4$}\,.
   \end{array} \right.
 $$
Each $U_-^{(z_i)}$ is equipped with a built-in morphism
  $\pi_-^{(i)}:U_-^{(z_i)}\rightarrow Y$
 in such a way that, when both are put together,
  they glue to give the resolution $\pi_-:Y_-\rightarrow Y$.

Consider the lifting
  $\{\pi_-^{(i)\,\prime}: U_-^{(z_i)} \rightarrow W\, |\,
                                                  i\,=\,1\,, 4\,\}$
  of the atlas
  $\{\, \pi_-^{(i)}:U_-^{(z_i)}\rightarrow Y\,|\, i\,=\,1\,, 4\,\}$
  of $Y_-$
 that is given by the lifting
  $\{\pi_-^{(i)\,\prime}:U_-^{(z_i)} \rightarrow W_{ut}\subset W\, |\,
                                                i\,=\,1\,, 4\,\}$
  defined by

  {\small
  $$
   \begin{array}{ccccc}
    \pi_-^{(1)\,\prime,\sharp}   & :
      & a_1\,,\, a_2\,,\, a_3\,,\, a_4\,,\,
        \varepsilon_1\,,\, \varepsilon_2\,,\,
        \varepsilon_3\,,\, \varepsilon_4
      & \longmapsto
      & z_1\,,\, z_3u_4\,,\, z_3\,,\, z_1u_4\,,\,
          1\,,\,      0\,,\,   0\,,\,    u_4
       \hspace{1em}\mbox{respectively}\,,\\[.6ex]
    \pi_-^{(4)\,\prime,\sharp}   & :
      & a_1\,,\, a_2\,,\, a_3\,,\, a_4\,,\,
        \varepsilon_1\,,\, \varepsilon_2\,,\,
        \varepsilon_3\,,\, \varepsilon_4
      & \longmapsto
      & z_4u_1\,,\, z_2\,,\, z_2u_1\,,\, z_4\,,\,
           u_1\,,\,   0\,,\,      0\,,\,   1
        \hspace{1em}\mbox{respectively}\,.
   \end{array}
  $$
  }

\noindent
The pair,  $\pi_-^{(1)\,\prime}$ and $\pi_-^{(4)\,\prime}$, are now
 embeddings into $W$ that are gluable
 up to the pointwise $\GL_2({\Bbb C})$-action.
Same construction as in the case of $\widetilde{Y}$ gives then
 a $\GL_2({\Bbb C})$-invariant subscheme $Y_-^{\prime}$ of $W$
 whose geometric quotient $Y_-^{\prime}/\GL_2({\Bbb C})$ exists
  and is equipped with a built-in isomorphism
  $Y_-^{\prime}/\GL_2({\Bbb C})\stackrel{\sim}{\rightarrow} Y_-$
 as schemes over $Y$.

This concludes the explicit construction.

\smallskip

\begin{sremark}
{\it $[\,$lifting to jet-scheme$\,]$.} {\rm
 Note that there is a one-to-one correspondence between
  $\GL_2({\Bbb C})$-orbits in $W$ and
  isomorphism classes of $0$-dimensional torsion sheaves
   of length $2$ on the conifold $Y$
   (i.e.\ the push-forward Chan-Paton sheaves on $Y$
    under associated morphisms from the Azumaya point
    $\Space M_2({\Bbb C})$ with the fundamental module ${\Bbb C}^2$)
  with connected support.
 Under this correspondence,
  the various special liftings-to-$W$ in the construction above:
   $$
    (\pi^{(1)\,\prime}\,,\, \pi^{(2)\,\prime}\,,\,
        \pi^{(3)\,\prime}\,,\, \pi^{(4)\,\prime})\,, \hspace{1em}
    (\pi_+^{(1)\,\prime}\,,\,\pi_+^{(3)\,\prime})\,, \hspace{1em}
    (\pi_-^{(1)\,\prime}\,,\,\pi_-^{(4)\,\prime})\,,
   $$
   and the gluing property, up to the pointwise $\GL_2({\Bbb C})$-action,
        in each tuple
  follow from the underlying lifting property to the related jet-schemes,
   which is the total space of the tangent sheaf ${\cal T}_Y$ of $Y$
   in our case.
}\end{sremark}

\bigskip

\begin{flushleft}
{\bf A comparison with resolutions via noncommutative desingularizations.}
\end{flushleft}
Consider the {\it conifold algebra} defined by\footnote{The
                          highlight here follows [leB-S]
                           with some change of notations for consistency
                           and mild rephrasings to link ibidem
                           directly with us.}
$$
 \Lambda_c\;:=\;
   \frac{{\Bbb C}\langle\,\xi_1\,,\, \xi_2\,,\, \xi_3\,\rangle}
        {(\,\xi_1^2\xi_2-\xi_2\xi_1^2\,,\, \xi_1\xi_2^2-\xi_2^2\xi_1\,,\,
            \xi_1\xi_3+\xi_3\xi_1\,,\,     \xi_2\xi_3+\xi_3\xi_2\,,\,
            \xi_3^2-1\,)}\,,
$$
where the numerator is the associative unital ${\Bbb C}$-algebra
 generated by $\{\xi_1,\, \xi_2,\, \xi_3\}$ and
 the denominator is the two-sided ideal generated by the elements
 of ${\Bbb C}\langle\xi_1\, \xi_2,\, \xi_3\rangle$ as indicated.

\smallskip

\begin{slemma}
{\bf [center of $\Lambda_c$].} {\rm ([leB-S: Lemma 5.4].)}
 The ${\Bbb C}$-algebra monomorphism
  $$
   \begin{array}{ccccc}
    \tau^{\sharp} & : & {\Bbb C}[z_1, z_2, z_3, z_4]/(z_1z_2-z_3z_4)
                  & \longrightarrow  & \Lambda_c   \\[.6ex]
    && z_1 & \longmapsto  & \xi_1^2                          \\[.6ex]
    && z_2 & \longmapsto  & \xi_2^2                          \\[.6ex]
    && z_3 & \longmapsto
     & \frac{1}{2}(\xi_1\xi_2+\xi_2\xi_1)\,
        +\,\frac{1}{2}(\xi_1\xi_2-\xi_2\xi_1)\xi_3 \\[.6ex]
    && z_4 & \longmapsto
     & \frac{1}{2}(\xi_1\xi_2+\xi_2\xi_1)\,
        -\,\frac{1}{2}(\xi_1\xi_2-\xi_2\xi_1)\xi_3
   \end{array}
  $$
  realizes ${\Bbb C}[z_1, z_2, z_3, z_4]/(z_1z_2-z_3z_4)$
  as the center of $\Lambda_c$.
\end{slemma}

\smallskip

\begin{sproposition}
{\bf [representation variety of $\Lambda_c$].}
{\rm ([leB-S: Proposition 5.7].)}
 The representation variety
  $\Rep(\Lambda_c, M_2({\Bbb C}))$
  is a smooth affine variety with three disjoint irreducible components.
 Two of these components are a point.
 The third $\Rep^0(\Lambda_c,M_2({\Bbb C}))$ has dimension $6$.
\end{sproposition}

\smallskip

This implies\footnote{Readers
                are referred to [leB1] for a general study of
                 the several notions involved in this paragraph.
                We do not need their details here.}
 that $\Lambda_c$ is a {\it smooth order} over
 ${\Bbb C}[z_1, z_2, z_3, z_4]/(z_1z_2-z_3z_4)$  and,
if one defines
  $\Spec \Lambda_c$ to be the set of two-sided prime ideals
   of $\Lambda_c$ with the Zariski topology,
 then the natural morphism
  $$
   \Spec \Lambda_c\; \longrightarrow\;
    \Spec({\Bbb C}[z_1,z_2,z_3,z_4]/(z_1z_2-z_3z_4))
  $$
  by intersecting a two-sided prime ideal of $\Lambda_c$
   with the center of $\Lambda_c$
  gives a {\it smooth noncommutative desingularization} of $Y$.
([leB-S: Proposition 5.7].)

Up to the conjugation by an element in $\GL_2({\Bbb C})$,
 a ${\Bbb C}$-algebra homomorphism
 $\rho:\Lambda_c\rightarrow M_2({\Bbb C})$
 can be put into one the following three forms:
 (In (1) and (2) below, $0$ and $\Id$ are respectively
  the zero matrix and the identity matrix in $M_2({\Bbb C})$.)
 \begin{itemize}
  \item[(1)]
   $\rho(\xi_1)=0\,$, $\rho(\xi_2)=0\,$,
    $\rho(\xi_3)=\hspace{1.4ex}\Id\,$;

  \item[(2)]
   $\rho(\xi_1)=0\,$, $\rho(\xi_2)=0\,$, $\rho(\xi_3)=-\Id\,$;

  \item[(3)]
   $$
    \rho(\xi_1)\; =\;
     \left[ \begin{array}{cc} 0 & a_1 \\ b_1 & 0 \end{array} \right]\,,
     \hspace{1em}
    \rho(\xi_2)\; =\;
     \left[ \begin{array}{cc} 0 & a_2 \\ b_2 & 0 \end{array} \right]\,,
     \hspace{1em}
    \rho(\xi_3)\; =\;
     \left[ \begin{array}{cc} 1 & 0 \\ 0 & -1 \end{array} \right]\,.
   $$
 \end{itemize}
 Form (1) and Form (2) correspond to the two point-components
  in $\Rep(\Lambda_c,M_2({\Bbb C}))$  and
 Form (3) corresponds to elements in $\Rep^0(\Lambda_c,M_2({\Bbb C}))$.
 On the subvariety ${\Bbb A}^4_{[a_1, b_1, a_2, b_2]}$
  of $\Rep^0(\Lambda_c, M_2({\Bbb C}))$ that parameterizes
  $\rho$ of the form (3),
 the $\GL_2({\Bbb C})$-action on $\Rep^0(\Lambda_c, M_2({\Bbb C}))$
  reduces to the ${\Bbb C}^{\ast}\times {\Bbb C}^{\ast}$-action
  $$
   \xymatrix{
   (a_1, b_1, a_2, b_2) \ar[rrr]^-{(t_1, t_2)}
    &&& (\, t_1t_2^{-1}a_1\,,\, t_1^{-1}t_2b_1\,,\,
          t_1t_2^{-1}a_2\,,\, t_1^{-1}t_2b_2\,)\,,
   }
  $$
  where $(t_1, t_2)\in {\Bbb C}^{\ast}\times{\Bbb C}^{\ast}$.
 The pair $(\rho(\xi_1), \rho(\xi_2))$ in Form (3) realizes
  this ${\Bbb A}^4_{[a_1, b_1, a_2, b_2]}$
  as the representation variety of the quiver
  $$
   \xymatrix{
    \bullet \ar @/_1ex/[rr]|{\;b_1} \ar @/_3ex/[rr]|{\;b_2}
     && \bullet \ar @/_3ex/[ll]|{\;a_1} \ar @/_1ex/[ll]|{\;a_2}\, .
   }
  $$

 Impose the trivial $\GL_2({\Bbb C})$-action on $Y$,
 then note that
  there is a natural $\GL_2({\Bbb C})$-equivariant morphism
  from $\Rep(\Lambda_c, M_2({\Bbb C}))$ to $Y$,
  as the composition
  $$
   {\Bbb C}[z_1,z_2, z_3, z_4]/(z_1z_2-z_3z_4)\;
    \stackrel{\tau^{\sharp}}{\longrightarrow}\; \Lambda_c\;
    \stackrel{\rho}{\longrightarrow}\; M_2({\Bbb C})
  $$
  has the form
  $$
    z_i\; \longmapsto\; 0\;,
          \hspace{1em}i\,=\, 1\,,\, 2\,,\, 3\,,\, 4\,,
  $$
  for $\rho$ conjugate to Form (1) or Form (2);
  $$
   z_1\;\longmapsto\; a_1b_1\,\Id\,,\hspace{1em}
   z_2\;\longmapsto\; a_2b_2\,\Id\,,\hspace{1em}
   z_3\;\longmapsto\; a_1b_2\,\Id\,,\hspace{1em}
   z_4\;\longmapsto\; a_2b_1\,\Id\,
  $$
  for $\rho$ conjugate to Form (3).\footnote{Note
                  that when restricted to
                   ${\Bbb A}^4_{[a_1, b_1, a_2, b_2]}
                    \subset${\it Rep}$\,^0(\Lambda_c, M_2({\Bbb C}))$,
                   this is the morphism
                   ${\Bbb A}^4_{[\xi_1, \xi_2,\xi_3,\xi_4]}\rightarrow Y$
                   in Sec.~2 after the substitution:
                    $a_1$ (here) $\rightarrow \xi_1$ (there),
                    $a_2 \rightarrow \xi_2$,
                    $b_1 \rightarrow \xi_3$,
                    $b_2 \rightarrow \xi_4$.}
 One can now follow the setting of [Ki] to define the stable structures
  for the $\GL_2({\Bbb C})$-action on $\Rep^0(\Lambda_c,M_2({\Bbb C}))$.
 There are two different choices,  $\theta_+$ and $\theta_-$,
  of such structures in the current case.
 The corresponding stable locus on the quiver variety
  ${\Bbb A}^4_{[a_1, b_1, a_2, b_2]}$
  is given respectively by
  $$
   {\Bbb A}_{[a_1, b_1, a_2, b_2]}^{4\,,\,\theta_+}\;
    =\;{\Bbb A}^4_{[a_1, b_1, a_2, b_2]}-V(b_1, b_2)
     \hspace{1em}\mbox{and}\hspace{1em}
   {\Bbb A}_{[a_1, b_1, a_2, b_2]}^{4\,,\,\theta_-}\;
    =\;{\Bbb A}^4_{[a_1, b_1, a_2, b_2]}-V(a_1, a_2)\,,
  $$
  where $V(a_1,a_2)$ (resp.\ $V(b_1, b_2)$)
   is the subvariety of ${\Bbb A}^4_{[a_1, b_1, a_2, b_2]}$
   associated to the ideal $(a_1, a_2)$ (resp.\ $(b_1, b_2)$).
 The corresponding GIT quotients
  $$
   \xymatrix{
    \Rep^0(\Lambda_c, M_2({\Bbb C}))/\!/^{\theta_+}\GL_2({\Bbb C})
     \ar[rd]_{\pi^{\theta_+}}
    && \Rep^0({\Lambda_c, M_2({\Bbb C})})
         /\!/^{\theta_-}\GL_2({\Bbb C})
          \ar[ld]^{\pi^{\theta_-}}                 \\
    & Y
   }
  $$
  recover
 $$
  \xymatrix{
   Y_+\ar[rd]_{\pi_+} && Y_-\ar[ld]^{\pi_-}\\
   & Y &&
   }
 $$
 at the beginning of the section.
See [leB-S], [leB2] for the mathematical detail and
 [Be], [B-L], [K-W] for the SQFT/stringy origin.

{From} the viewpoint of the Polchinski-Grothendieck Ansatz,
{\it both} the Azumaya-type noncommutative structure on D-branes and
 a noncommutative structure over $Y$ described by
 $\Space\Lambda_c$ come into play in the above setting.
As indicated by the explicit expression for
  $\rho\circ\tau^{\sharp}$ above,
 any morphism
  $\tilde{\varphi}: \Space M_2({\Bbb C})\rightarrow \Space\Lambda_c$
  has the property:
 \begin{itemize}
  \item[$\cdot$]
   The composition
   $$
    \Space M_2({\Bbb C})\;
     \stackrel{\tilde{\varphi}}{\longrightarrow}\; \Space\Lambda_c \;
     \stackrel{\tau}{\longrightarrow}\;  Y
   $$
   is a morphism $\varphi:=\tilde{\varphi}\circ\tau$
   from the Azumaya point $\pt^{A\!z}=\Space M_2({\Bbb C})$ to $Y$
   with the associated surrogate $\pt_{\varphi}\simeq \Spec{\Bbb C}$.
 \end{itemize}
Thus, the new ingredient of target-space noncommutativity
 comes into play as another key role toward resolutions of $Y$
 in the above setting
while the generalized-jet-resolution-of-singularity picture
 in our earlier discussion disappears.

\smallskip

\begin{sremark}
{\it $[\,$world-volume noncommutativity
          vs.\ target-space(-time) noncommutativity$\,]$.}
{\rm
 Such a ``trading"
   between a noncommutativity target and
    morphisms from Azumaya schemes to a commutative target
  suggests a partial duality
   between D-brane world-volume noncommutativity and
   target space(-time) noncommutativity.
} \end{sremark}

\smallskip

\noindent
{\sc Figure}~3-1.

\begin{figure}[htbp]
 \epsfig{figure=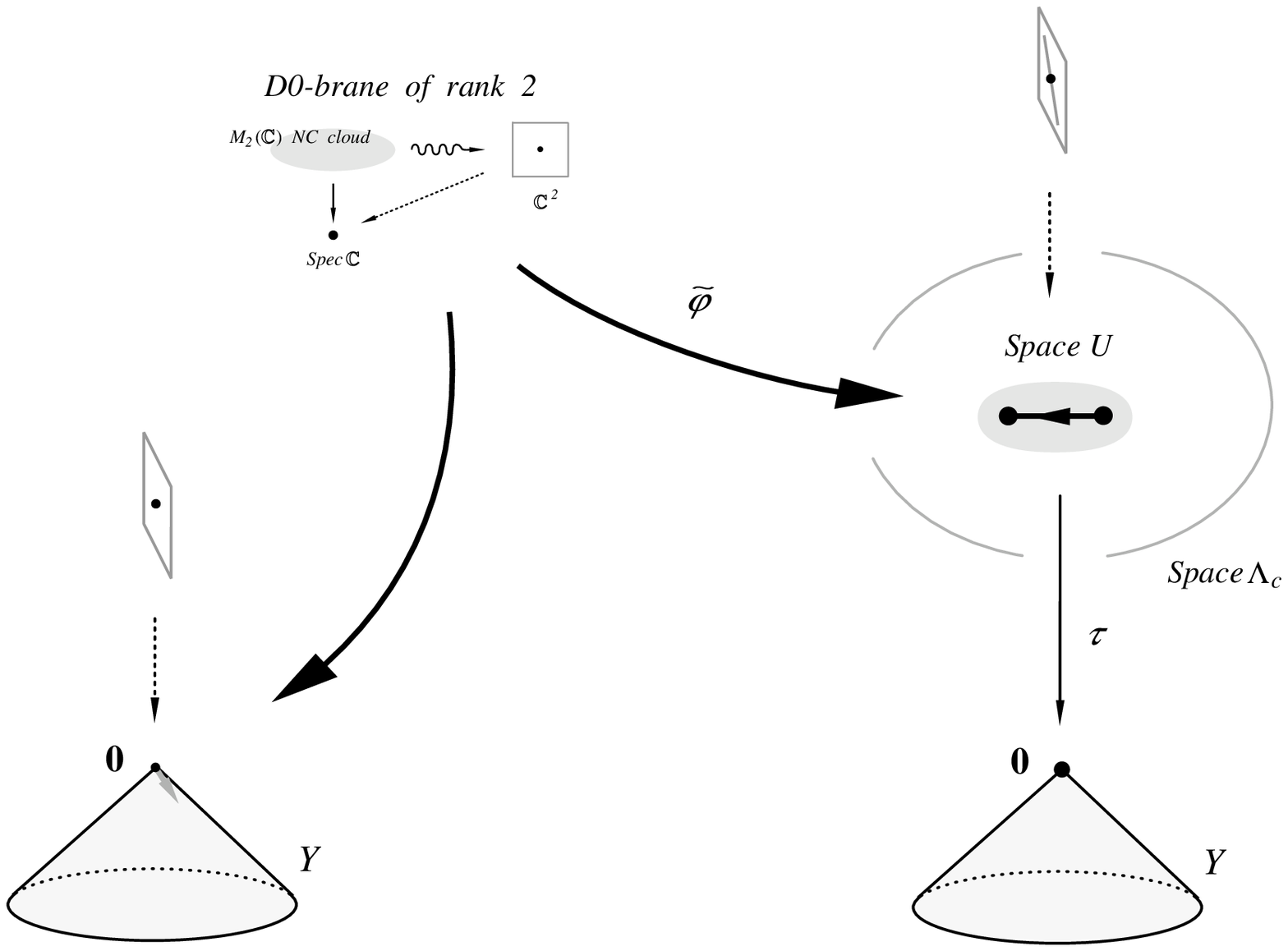,width=16cm}
 \centerline{\parbox{13cm}{\small\baselineskip 12pt
  {\sc Figure}~3-1.
 Trading of morphisms from {\it Space}$\,M_2({\Bbb C})$
   directly to the conifold $Y$
  with those to the noncommutative space {\it Space}$\,\Lambda_c$
   over $Y$.
 Note that for generic $\rho\in$ {\it Rep}$\,(\Lambda_c, M_2({\Bbb C}))$
  such that $\rho\circ \tau^{\sharp}=0$,
  $\rho(\Lambda_c)$ is similar to
  the ${\Bbb C}$-subalgebra $U$ of upper triangular matrices in
  $M_2({\Bbb C})$.
 The noncommutative point {\it Space}$\,U$ is also smooth,
  with {\it Spec}$\,U$ consisting of two ${\Bbb C}$-points
  connected by a directed nilpotent bond.
 It is thus represented by a quiver
  $\,\xymatrix{\bullet \ar[r] &\bullet}\,$ in the figure.
 Furthermore, let
  $\widetilde{\varphi}:$ {\it Space}$\,M_2({\Bbb C})\rightarrow$
   {\it Space}$\,\Lambda_c$
  be the corresponding morphism.
 Then $\widetilde{\varphi}$ determines also a flag
  in the Chan-Paton module $\widetilde{\varphi}_{\ast}{\Bbb C}^2$
  on the image D$0$-brane {\it Im}$\,\widetilde{\varphi}\,$.
 On the other hand, over a generic $p\ne {\mathbf 0}$ on $Y$,
  the generic image of a $\widetilde{\varphi}^{\prime}$
  that maps to $p$ after the composition with $\tau$
  will be simply {\it Space}$\,M_2({\Bbb C})\,$.
 }}
\end{figure}

\bigskip

\newpage
\baselineskip 13pt
%references
%endfootnotesize

\end{document}

%% file: newcmd.tex
\newcommand{\AffineSpace}{\mbox{\it ${\cal A}$ffine${\cal S}\!$pace}\,}
\newcommand{\AlgCategory}{\mbox{\it ${\cal A}$lg}\,}

 \newcommand{\Algscriptsize}{\mbox{\scriptsize\it Alg}\,}

\newcommand{\Bl}{\mbox{\it Bl}\,}

\newcommand{\End}{\mbox{\it End}\,}

\newcommand{\GL}{\mbox{\it GL}}

\newcommand{\Hom}{\mbox{\it Hom}\,}

\newcommand{\Id}{\mbox{\it Id}}

\newcommand{\Image}{\mbox{\it Im}\,}

\newcommand{\Ker}{\mbox{\it Ker}\,}

\newcommand{\Mor}{\mbox{\it Mor}\,}

\newcommand{\NA}{\mbox{\it N${\Bbb A}$}}

\newcommand{\Proj}{\mbox{\it Proj}\,}

\newcommand{\PGL}{\mbox{\it PGL}\,}

\newcommand{\Rep}{\mbox{\it Rep}\,}
 \newcommand{\scriptsizeRep}{\mbox{\scriptsize\it Rep}\,}

\newcommand{\SetCategory}{\mbox{\it ${\cal S}\!$et}\,}

\newcommand{\Space}{\mbox{\it Space}\,}

\newcommand{\Span}{\mbox{\it Span}\,}
\newcommand{\Spec}{\mbox{\it Spec}\,}

\newcommand{\SU}{\mbox{\it SU}}

\newcommand{\Sym}{\mbox{\it Sym}}

\newcommand{\Tr}{\mbox{\it Tr}\,}

\newcommand{\pt}{\mbox{\it pt}}

\newcommand{\redscriptsize}{\mbox{\scriptsize\rm red}\,}

\newcommand{\singletonscriptsize}{\mbox{\scriptsize\it singleton}\,}